\documentclass[11pt]{article}%
\usepackage{amssymb}
\usepackage{amsmath}
\usepackage{amsthm}
\usepackage{graphicx}
\usepackage{amsfonts}%
\providecommand{\U}[1]{\protect\rule{.1in}{.1in}}
\theoremstyle{plain}
\newtheorem*{GL}{The grand plan cannot work for $p=3$}
\newtheorem*{SGT}{Sophie Germain's Theorem}
\newtheorem*{LS}{Theorem (Large Size of Solutions)}
\newtheorem*{NC}{Condition N-C}
\newtheorem*{2Np}{Condition $2$-N-$p$}
\newtheorem*{pNp}{Condition $p$-N-$p$}
\newtheorem*{Np-inv}{Condition N-$p^{-1}$}
\theoremstyle{definition}
\newtheorem*{pf}{Proof}
\theoremstyle{remark}

\begin{document}

\title{\textquotedblleft Voici ce que j'ai trouv\'{e}:\textquotedblright\ Sophie
Germain's grand plan to prove Fermat's Last Theorem}
\author{Reinhard Laubenbacher\\Virginia Bioinformatics Institute\\Virginia Polytechnic Institute and State University\\Blacksburg, VA 24061, USA\\reinhard@vbi.vt.edu
\and David Pengelley\thanks{Dedicated to the memory of my parents, Daphne and Ted
Pengelley, for inspiring my interest in history, and to Pat Penfield, for her
talented, dedicated, and invaluable editorial help, love and enthusiasm, and
support for this project.}\\Mathematical Sciences\\New Mexico State University \\Las Cruces, NM 88003, USA\\davidp@nmsu.edu}
\maketitle

\begin{abstract}
A study of Sophie Germain's extensive manuscripts on Fermat's Last Theorem
calls for a reassessment of her work in number theory. There is much in these
manuscripts beyond the single theorem for Case 1 for which she is known from a
published footnote by Legendre. Germain had a fully-fledged, highly developed,
sophisticated plan of attack on Fermat's Last Theorem. The supporting
algorithms she invented for this plan are based on ideas and results
discovered independently only much later by others, and her methods are quite
different from any of Legendre's. In addition to her program for proving
Fermat's Last Theorem in its entirety, Germain also made major efforts at
proofs for particular families of exponents. The isolation Germain worked in,
due in substantial part to her difficult position as a woman, was perhaps
sufficient that much of this extensive and impressive work may never have been
studied and understood by anyone.\medskip

Une \'{e}tude approfondie des manuscrits de Sophie Germain sur le dernier
th\'{e}or\`{e}me de Fermat, r\'{e}v\`{e}le que l'on doit r\'{e}\'{e}valuer ses
travaux en th\'{e}orie des nombres. En effet, on trouve dans ses manuscrits
beaucoup plus que le simple th\'{e}or\`{e}me du premier cas que Legendre lui
avait attribu\'{e} dans une note au bas d'une page et pour lequel elle est
reconnue. Mme Germain avait un plan tr\`{e}s \'{e}labor\'{e} et
sophistiqu\'{e} pour prouver enti\`{e}rement ce dernier th\'{e}or\`{e}me de
Fermat. Les algorithmes qu'elle a invent\'{e}s sont bas\'{e}s sur des
id\'{e}es et resultats qui ne furent ind\'{e}pendamment d\'{e}couverts que
beaucoup plus tard. Ses m\'{e}thodes sont \'{e}galement assez diff\'{e}rentes
de celles de Legendre. En plus, Mme Germain avait fait de remarquables
progr\`{e}s dans sa recherche concernant certaines familles d'exposants.
L'isolement dans lequel Sophie Germain se trouvait, en grande partie d\^{u} au
fait qu'elle \'{e}tait une femme, fut peut-\^{e}tre suffisant, que ses
impressionnants travaux auraient pu passer compl\`{e}tement inaper\c{c}us et
demeurer incompris.

\medskip

\noindent\emph{MSC:} 01A50; 01A55; 11-03; 11D41

\medskip

\noindent\emph{Keywords: }Sophie Germain; Fermat's Last Theorem; Adrien-Marie
Legendre; Carl Friedrich Gauss; Guglielmo (Guillaume) Libri; number theory

\end{abstract}
\tableofcontents

\section{Introduction}

\begin{figure}
[pt]
\begin{center}
\includegraphics[
natheight=3.335600in,
natwidth=2.277900in,
height=3.3823in,
width=2.3177in
]%
{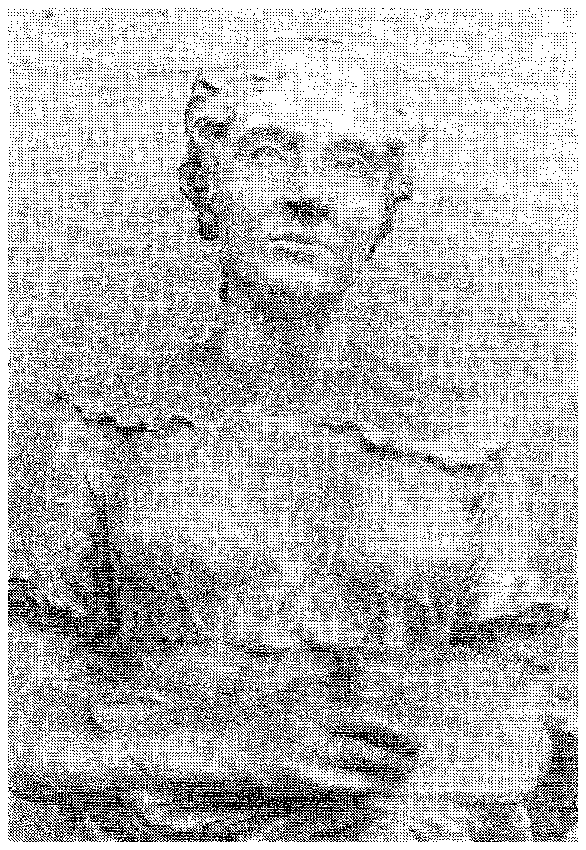}%
\caption{Sophie Germain: a bust by Z. Astruc}%
\label{germain-600-ph-image}%
\end{center}
\end{figure}

Sophie Germain (Figure \ref{germain-600-ph-image})\footnote{From \cite[p.
17]{bucc}.} was the first woman known for important original research in
mathematics.\footnote{A biography of Germain, with concentration on her work
in elasticity theory, discussion of her personal and professional life, and
references to the historical literature about her, is the book by Lawrence
Bucciarelli and Nancy Dworsky \cite{bucc}.} While perhaps more famous for her
work in mathematical physics that earned her a French Academy prize, Germain
is also credited with an important result in number theory towards proving
Fermat's Last Theorem. We will make a substantial reevaluation of her work on
the Fermat problem, based on translation and mathematical interpretation of
numerous documents in her own hand, and will argue that her accomplishments
are much broader, deeper, and more significant than has been realized.

Fermat's Last Theorem refers to Pierre de Fermat's famous seventeenth century
claim that the equation $z^{p}=x^{p}+y^{p}$ has no natural number solutions
$x$, $y$, $z$ for natural number exponents $p>2$. The challenge of proving
this assertion has had a tumultuous history, culminating in Andrew Wiles'
success at the end of the twentieth century \cite[XI.2]{ribenboim-amateurs}.

Once Fermat had proven his claim for exponent $4$ \cite[p. 615ff]{dickson}
\cite[p. 75ff]{weil}, it could be fully confirmed just by substantiating it
for odd prime exponents. But when Sophie Germain began working on the problem
at the turn of the nineteenth century, the only prime exponent that had a
proof was $3$ \cite[XXVI]{dickson} \cite[ch. 3]{edwards} \cite[I.6,
IV]{ribenboim-amateurs} \cite[p. 335ff]{weil}. As we will see, Germain not
only developed the one theorem she has long been known for towards proving
part of Fermat's Last Theorem for all primes. Her manuscripts reveal a
comprehensive program to prove Fermat's Last Theorem in its entirety.

\subsection{Germain's background and mathematical development}%

Sophie Germain\footnote{Much of our description here of Germain's background
appears also in \cite{musielak-review}.} was born on April 1, 1776 and lived
with her parents and sisters in the center of Paris throughout the upheavals
of the French Revolution. Even if kept largely indoors, she must as a teenager
have heard, and perhaps seen, some of its most dramatic and violent events.
Moreover, her father, Ambroise-Fran\c{c}ois Germain, a silk merchant, was an
elected member of the third estate to the Constituent Assembly convened in
1789 \cite[p. 9ff]{bucc}. He thus brought home daily intimate knowledge of
events in the streets, the courts, etc.; how this was actually shared, feared,
and coped with by Sophie Germain and her family we do not know.

Much of what we know of Germain's life comes from the biographical obituary
\cite{libri-notice} published by her friend and fellow mathematician Guglielmo
Libri, shortly after her death in 1831. He wrote that at age thirteen, Sophie
Germain, partly as sustained diversion from her fears of the Revolution
beginning outside her door, studied first Montucla's \emph{Histoire des
math\'{e}matiques}, where she read of the death of Archimedes on the sword of
a Roman soldier during the fall of Syracuse, because he could not be
distracted from his mathematical meditations. And so it seems that Sophie
herself followed Archimedes, becoming utterly absorbed in learning
mathematics, studying without any teacher from a then common mathematical work
by \'{E}tienne Bezout that she found in her father's library.

Her family at first endeavored to thwart her in a taste so unusual and
socially unacceptable for her age and sex. According to Libri, Germain rose at
night to work from the glimmer of a lamp, wrapped in covers, in cold that
often froze the ink in its well, even after her family, in order to force her
back to bed, had removed the fire, clothes, and candles from her room; it is
thus that she gave evidence of a passion that they thereafter had the wisdom
not to oppose further. Libri writes that one often heard of the happiness with
which Germain rejoiced when, after long effort, she could persuade herself
that she understood the language of analysis in Bezout. Libri continues that
after Bezout, Germain studied Cousin's differential calculus, and was absorbed
in it during the Reign of Terror in 1793--1794. It is from roughly 1794
onwards that we have some records of Germain interacting with the public
world. It was then, Libri explains, that Germain did something so rashly
remarkable that it would actually lack believability if it were mere fiction.

Germain, then eighteen years old, first somehow obtained the lesson books of
various professors at the newly founded \'{E}cole Polytechnique, and was
particularly focused on those of Joseph-Louis Lagrange on analysis. The
\'{E}cole, a direct outgrowth of the French Revolution, did not admit women,
so Germain had no access to this splendid new institution and its faculty.
However, the \'{E}cole did have the novel feature, heralding a modern
university, that its professors were both teachers and active researchers.
Indeed its professors included some of the best scientists and mathematicians
in the world. Libri writes that professors had the custom, at the end of their
lecture courses, of inviting their students to present them with written
observations. Sophie Germain, assuming the name of an actual student at the
\'{E}cole Polytechnique, one Antoine-August LeBlanc, submitted her
observations to Lagrange, who praised them, and learning the true name of the
imposter, actually went to her to attest his astonishment in the most
flattering terms.

Perhaps the most astounding aspect is that Germain appears to have entirely
self-educated herself to at least the undergraduate level, capable of
submitting written student work to Lagrange, one of the foremost researchers
in the world, that was sufficiently notable to make him seek out the author.
Unlike other female mathematicians before her, like Hypatia, Maria Agnesi, and
\'{E}milie du Ch\^{a}telet, who had either professional mentors or formal
education to this level, Sophie Germain appears to have climbed to university
level unaided and entirely on her own initiative.

Libri continues that Germain's appearance thus on the Parisian mathematical
scene drew other scholars into conversation with her, and that she became a
passionate student of number theory with the appearance of Adrien-Marie
Legendre's (Figure \ref{legendre-300-ph-image}) \emph{Th\'{e}orie des Nombres}
in 1798. Both Lagrange and Legendre became important personal mentors to
Germain, even though she could never attend formal courses of study. After
Carl Friedrich Gauss's \emph{Disquisitiones Arithmeticae }appeared in 1801,
Germain took the additional audacious step, in 1804, of writing to him, again
under the male pseudonym of LeBlanc (who in the meantime had died), enclosing
some research of her own on number theory, and particularly on Fermat's Last
Theorem. Gauss entered into serious mathematical correspondence with
\textquotedblleft Monsieur LeBlanc\textquotedblright. In 1807 the true
identity of LeBlanc was revealed to Gauss when Germain intervened with a
French general, a family friend, to ensure Gauss's personal safety in
Braunschweig during Napoleon's Jena campaign \cite[ch. 2, 3]{bucc}. Gauss's
response to this surprise metamorphosis of his correspondent was
extraordinarily complimentary and encouraging to Germain as a mathematician,
and quite in contrast to the attitude of many 19th century scientists and
mathematicians about women's abilities:%

\begin{figure}
[pt]
\begin{center}
\includegraphics[
natheight=3.451500in,
natwidth=2.936900in,
height=2.8582in,
width=2.437in
]%
{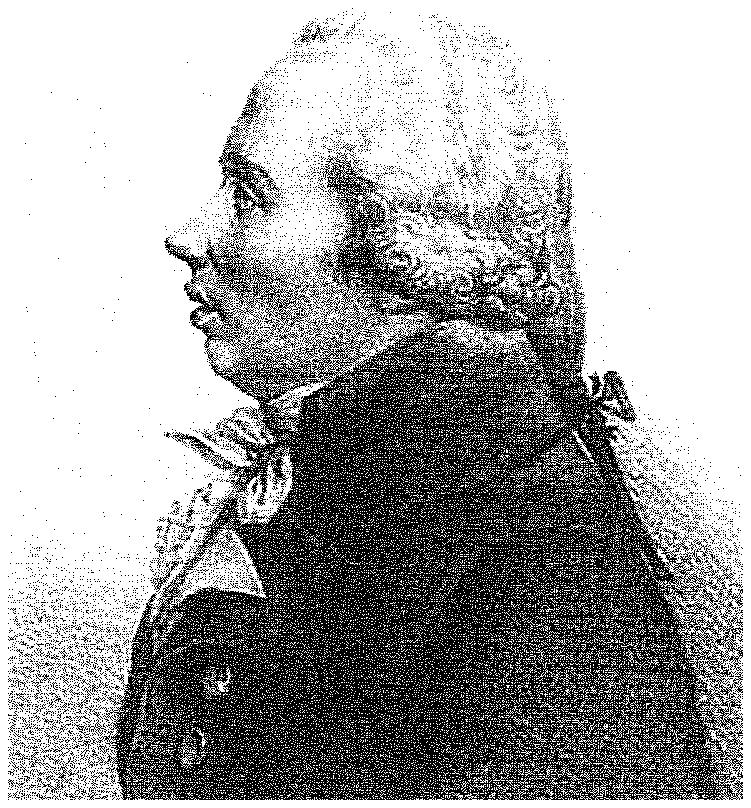}%
\caption{Adrien-Marie Legendre}%
\label{legendre-300-ph-image}%
\end{center}
\end{figure}

\begin{quote}
But how can I describe my astonishment and admiration on seeing my esteemed
correspondent Monsieur LeBlanc metamorphosed into this celebrated person,
yielding a copy so brilliant it is hard to believe? The taste for the abstract
sciences in general and, above all, for the mysteries of numbers, is very
rare: this is not surprising, since the charms of this sublime science in all
their beauty reveal themselves only to those who have the courage to fathom
them. But when a woman, because of her sex, our customs and prejudices,
encounters infinitely more obstacles than men, in familiarizing herself with
their knotty problems, yet overcomes these fetters and penetrates that which
is most hidden, she doubtless has the most noble courage, extraordinary
talent, and superior genius. Nothing could prove to me in a more flattering
and less equivocal way that the attractions of that science, which have added
so much joy to my life, are not chimerical, than the favor with which you have
honored it.

The scientific notes with which your letters are so richly filled have given
me a thousand pleasures. I have studied them with attention and I admire the
ease with which you penetrate all branches of arithmetic, and the wisdom with
which you generalize and perfect. \cite[p.\ 25]{bucc}
\end{quote}

The subsequent arcs of Sophie Germain's two main mathematical research
trajectories, her interactions with other researchers, and with the
professional institutions that forced her, as a woman, to remain at or beyond
their periphery, are complex. Germain's development of a mathematical theory
explaining the vibration of elastic membranes is told by Bucciarelli and
Dworsky in their mathematical biography \cite{bucc}. In brief, the German
physicist Ernst Chladni created a sensation in Paris in 1808 with his
demonstrations of the intricate vibrational patterns of thin plates, and at
the instigation of Napoleon, the Academy of Sciences set a special prize
competition to find a mathematical explanation. Germain pursued a theory of
vibrations of elastic membranes, and based on her partially correct
submissions, the Academy twice extended the competition, finally awarding her
the prize in 1816, while still criticizing her solution as incomplete, and did
not publish her work \cite[ch. 7]{bucc}. The whole experience was definitely
bittersweet for Germain.

The Academy then immediately established a new prize, for a proof of Fermat's
Last Theorem. While Sophie Germain never submitted a solution to this new
Academy prize competition, and she never published on Fermat's Last Theorem,
it has long been known that she worked on it, from the credit given her in
Legendre's own 1823 memoir published on the topic \cite[p. 87]{bucc} \cite[p.
189]{lp} \cite{legendre}. Our aim in this paper is to analyze the surprises
revealed by Germain's manuscripts and letters, containing work on Fermat's
Last Theorem going far beyond what Legendre implies.

We will find that the results Legendre credits to Germain were merely a small
piece of a much larger body of work. Germain pursued nothing less than an
ambitious full-fledged plan of attack to prove Fermat's Last Theorem in its
entirety, with extensive theoretical techniques, side results, and supporting
algorithms. What Legendre credited to her, known today as Sophie Germain's
Theorem, was simply a small part of her big program, a piece that could be
encapsulated and applied separately as an independent theorem, as was put in
print by Legendre.

\subsection{Germain's number theory in the literature}

Sophie Germain's principal work on the Fermat problem has long been understood
to be entirely described by a single footnote in Legendre's 1823 memoir
\cite[p. 734]{dickson} \cite[ch. 3]{edwards} \cite[\S 22]{legendre} \cite[p.
110]{ribenboim-amateurs}. The memoir ends with Legendre's own proof for
exponent $5$, only the second odd exponent for which it was proven. What
interests us here, though, is the first part of his treatise, since Legendre
presents a general analysis of the Fermat equation whose main theoretical
highlight is a theorem encompassing all odd prime exponents, today named after Germain:

\begin{SGT}
For an odd prime exponent $p$, if there exists an auxiliary prime $\theta$
such that there are no two nonzero consecutive $p$-th powers modulo $\theta$,
nor is $p$ itself a $p$-th power modulo $\theta$, then in any solution to the
Fermat equation $z^{p}=x^{p}+y^{p}$, one of $x$, $y$, or $z$ must be divisible
by $p^{2}$.
\end{SGT}

Sophie Germain's Theorem can be applied for many prime exponents, by producing
a valid auxiliary prime, to eliminate the existence of solutions to the Fermat
equation involving numbers not divisible by the exponent $p$. This elimination
is today called Case 1\ of Fermat's Last Theorem. Work on Case 1 has continued
to the present, and major results, including for instance its recent
establishment for infinitely many prime exponents \cite{adleman,fouvry}, have
been proven by building on the very theorem that Germain introduced.

Before proceeding further, we briefly give the minimum mathematical background
needed to understand fully the statement of the theorem, and then an
illustration of its application. The reader familiar with modular arithmetic
may skip the next two paragraphs.

Two whole numbers $a$ and $b$ are called \textquotedblleft
congruent\textquotedblright\ (or \textquotedblleft
equivalent\textquotedblright) \textquotedblleft modulo $\theta$%
\textquotedblright\ (where $\theta$ is a natural number called the modulus) if
their difference $a-b$ is a multiple of $\theta$; this is easily seen to
happen precisely if they have the same remainder (\textquotedblleft
residue\textquotedblright) upon division by $\theta$. (Of course the residues
are numbers between $0$ and $\theta-1$, inclusive.) We write $a\equiv b$
$\left(  \operatorname{mod}\theta\right)  $ and say \textquotedblleft$a$ is
congruent to $b$ modulo $\theta$\textquotedblright\ (or for short, just
\textquotedblleft$a$ is $b$ modulo $\theta$\textquotedblright).\footnote{The
notation and language of congruences was introduced by Gauss in his
\emph{Disquisitiones Arithmeticae} in 1801, and Sophie Germain was one of the
very first to wholeheartedly and profitably adopt it in her research.}
Congruence satisfies many of the same simple properties that equality of
numbers does, especially in the realms of addition, subtraction, and
multiplication, making it both useful and easy to work with. The reader will
need to become familiar with these properties, and we will not spell them out
here. The resulting realm of calculation is called \textquotedblleft modular
arithmetic\textquotedblright, and its interesting features depend very
strongly on the modulus $\theta$.

In the statement of the theorem, when one considers whether two numbers are
\textquotedblleft consecutive modulo $\theta$\textquotedblright, one means
therefore not that their difference is precisely $1$, but rather that it is
congruent to $1$ modulo $\theta$; notice that one can determine this by
looking at the residues of the two numbers and seeing if the residues are
consecutive. (Technically, one also needs to recognize as consecutive modulo
$\theta$ two numbers whose residues are $0$ and $\theta-1$, since although the
residues are not consecutive as numbers, the original numbers will have a
difference congruent to $0-\left(  \theta-1\right)  =1-\theta$, and therefore
to $1$ $\left(  \operatorname{mod}\theta\right)  $. In other words, the
residues $0$ and $\theta-1$ should be thought of as consecutive in how they
represent numbers via congruence. However, since we are interested only in
numbers with nonzero residues, this complication will not arise for us.)

We are ready for an example. Let us choose $p=3$ and $\theta=13$, both prime,
and test the two hypotheses of Sophie Germain's Theorem by brute force
calculation. We need to find all the nonzero residues of $3$rd powers (cubic
residues) modulo $13$. A basic feature of modular arithmetic tells us that we
need only consider the cubes of the possible residues modulo $13$, i.e., from
$0$ to $12$, since all other numbers will simply provide cyclic repetition of
what these produce. And since we only want nonzero results modulo $\theta$, we
may omit $0$. Brute force calculation produces Table \ref{cubicresidues}.%

\begin{table}[t]
\[
\begin{tabular}
[c]{|c|r|r|r|r|r|r|r|r|r|r|r|r|}\hline
Residue & $1$ & $2$ & $3$ & $4$ & $5$ & $6$ & $7$ & $8$ & $9$ & $10$ & $11$ &
$12$\\\hline
Cube & $1$ & $8$ & $27$ & $64$ & $125$ & $216$ & $343$ & $512$ & $729$ &
$1000$ & $1331$ & $1728$\\\hline\begin{tabular}
[c]{c}
Cubic\\
residue
\end{tabular}
& $1$ & $8$ & $1$ & $12$ & $8$ & $8$ & $5$ & $5$ & $1$ & $12$ & $5$ &
$12$\\\hline\end{tabular}
\]
\caption{Cubic residues modulo $13$}
\label{cubicresidues}
\end{table}%

For instance, the residue of $8^{3}=512$ modulo $13$ can be obtained by
dividing $512$ by $13$, with a remainder of $5$. However, there are much
quicker ways to obtain this, since in a congruence calculation, any number
(except exponents) may be replaced with anything congruent to it. So for
instance we can easily calculate that $8^{3}=64\cdot8\equiv\left(  -1\right)
\cdot\left(  -5\right)  =5$ $\left(  \operatorname{mod}13\right)  $.

Now we ask whether the two hypotheses of Sophie Germain's Theorem are
satisfied? Indeed, no pair of the nonzero cubic residues $1,$ $5,$ $8,$ $12$
modulo $13$ are consecutive, and $p=3$ is not itself among the residues. So
Sophie Germain's Theorem proves that any solution to the Fermat equation
$z^{3}=x^{3}+y^{3}$ would have to have one of $x$, $y$, or $z$ divisible by
$p^{2}=9$.

Returning to Legendre's treatise, after the theorem he supplies a table
verifying the hypotheses of the theorem for $p<100$ by brute force display of
all the $p$-th power residues modulo a single auxiliary prime $\theta$ chosen
for each value of $p$. Legendre then credits Sophie Germain with both the
theorem, which is the first general result about arbitrary exponents for
Fermat's Last Theorem, and its successful application for $p<100$. One assumes
from Legendre that Germain developed the brute force table of residues as her
means of verification and application of her theorem. Legendre continues on to
develop more theoretical means of verifying the hypotheses of Sophie Germain's
Theorem, and he also pushes the analysis further to demonstrate that any
solutions to the Fermat equation for certain exponents would have to be
extremely large.

For almost two centuries, it has been assumed that this theorem and its
application to exponents less than $100$, the basis of Germain's reputation,
constitute her entire contribution to Fermat's Last Theorem. However, we will
find that this presumption is dramatically off the mark as we study Germain's
letters and manuscripts. The reward is a wealth of new material, a vast
expansion over the very limited information known just from Legendre's
footnote. We will explore its enlarged scope and extent. Figures
\ref{figure-1} and \ref{figure-2} in Section \ref{precis} show all the
interconnected pieces of her work, and the place of Sophie Germain's Theorem
in it. The ambitiousness and importance of Germain's work will prompt a major
reevaluation, and recommend a substantial elevation of her reputation.

Before considering Germain's own writing, we note that the historical record
based solely on Legendre's footnote has itself been unjustly portrayed. Even
the limited results that Legendre clearly attributed to Germain have been
understated and misattributed in much of the vast secondary literature. Some
writers state only weaker forms of Sophie Germain's Theorem, such as merely
for $p=5$, or only for auxiliary primes of the form $2p+1$ (known as
\textquotedblleft Germain primes\textquotedblright, which happen always to
satisfy the two required hypotheses). Others only conclude divisibility by the
first power of $p$, and some writers have even attributed the fuller $p^{2}%
$-divisibility, or the determination of qualifying auxiliaries for $p<100$, to
Legendre rather than to Germain. A few have even confused the results Legendre
credited to Germain with a completely different claim she had made in her
first letter to Gauss, in 1804 \cite{stupuy}. We will not list all these
failings here. Fortunately a few books have correctly stated Legendre's
attribution to Germain \cite[p. 734]{dickson} \cite[ch. 3]{edwards} \cite[p.
110]{ribenboim-amateurs}. We will not elaborate in detail on the huge related
mathematical literature except for specific relevant comparisons of
mathematical content with Germain's own work. Ribenboim's most recent book
\cite{ribenboim-amateurs} gives a good overall history of related
developments, including windows into the intervening literature.

\subsection{Manuscript sources, recent research, and scope}

Bucciarelli and Dworsky's mathematical biography of Germain's work on
elasticity theory \cite{bucc} utilized numerous Germain manuscripts from the
archives of the Biblioth\`{e}que Nationale in Paris. Many other Germain
manuscripts are also held in the Biblioteca Moreniana in Firenze (Florence)
\cite[pp. 229--235, 239--241]{delcentinabook} \cite{delcentina-germain}%
.\footnote{The story of how Germain's manuscripts ended up in the two archives
is an extraordinary one, a consequence of the amazing career of Guglielmo
(Guillaume)\ Libri, mathematician, historian, bibliophile, thief, and friend
of Sophie Germain \cite{delcentinabook,ruju}.
\par
Exactly how Libri originally obtained Germain's manuscripts remains uncertain.
We note, however, that Germain was not affiliated with any institution that
might naturally have taken them, while Libri was a good friend of hers. After
his expulsion from Tuscany for his role in the plot to persuade the Grand-Duke
to promulgate a constitution, Libri traveled for many months, not reaching
Paris until fully six months after Germain died. Nonetheless, it seems he
ended up with almost all her papers \cite[p.~142f]{delcentinabook}, and it was
entirely in character for him to manage this, since he built a gargantuan
private library of important books, manuscripts, and letters
\cite{delcentinabook}.
\par
It appears that many of Germain's manuscripts in the Biblioth\`{e}que
Nationale were probably among those confiscated by the police from Libri's
apartment at the Sorbonne when he fled to London in 1848 to escape the charge
of thefts from French public libraries \cite[p. 146]{delcentinabook}. The
Germain manuscripts in the Biblioteca Moreniana were among those shipped with
Libri's still remaining vast collection of books and manuscripts before he set
out to return from London to Florence in 1868. These latter Germain materials
are among those fortunate to have survived intact despite a long and tragic
string of events following Libri's death in 1869
\cite{delcentinabook,delcentina-germain}. Ultimately it seems that Libri was
the good fortune that saved Germain's manuscripts; otherwise they might simply
have drifted into oblivion. See also
\cite{delcentina-abel1,delcentina-abel2,delcentina-fate} for the story of Abel
manuscripts discovered in the Libri collections in the Biblioteca Moreniana.}
While Bucciarelli and Dworsky focused primarily on her work on elasticity
theory, many of the manuscripts in these archives are on number theory. Their
book also indicates there are unpublished letters from Germain to Gauss, held
in G\"{o}ttingen; in particular, there is a letter written in 1819 almost
entirely about Fermat's Last Theorem.

It appears that Germain's number theory manuscripts have received little
attention during the nearly two centuries since she wrote them. We began
working with them in 1994, and published a translation and analysis of
excerpts from one (Manuscript B below) in our 1999 book \cite[p. 190f]{lp}. We
demonstrated there that the content and proof of Sophie Germain's Theorem, as
attributed to her by Legendre, is implicit within the much broader aims of
that manuscript, thus substantiating in Germain's own writings Legendre's
attribution. Since then we have analyzed the much larger corpus of her number
theory manuscripts, and we present here our overall evaluation of her work on
Fermat's Last Theorem, which forms a coherent theory stretching over several
manuscripts and letters.

Quite recently, and independently from us, Andrea Del Centina
\cite{delcentina-germain-flt} has also transcribed and analyzed some of
Germain's manuscripts, in particular one at the Biblioteca Moreniana and its
more polished copy at the Biblioth\`{e}que Nationale (Manuscripts D and A
below). While there is some overlap between Del Centina's focus and ours,
there are major differences in which manuscripts we consider, and in what
aspects of them we concentrate on. In fact our research and Del Centina's are
rather complementary in what they analyze and present. Overall there is no
disagreement between the main conclusions we and Del Centina draw; instead
they supplement each other. After we list our manuscript sources below, we
will compare and contrast Del Centina's specific selection of manuscripts and
emphasis with ours, and throughout the paper we will annotate any specifically
notable comparisons of analyses in footnotes.

Germain's handwritten papers on number theory in the Biblioth\`{e}que
Nationale are almost all undated, relatively unorganized, and unnumbered
except by the archive. And they range all the way from scratch paper to some
beautifully polished finished pieces. We cannot possibly provide a definitive
evaluation here of this entire treasure trove, nor of all the manuscripts in
the Biblioteca Moreniana. We will focus our attention within these two sets of
manuscripts on the major claims about Fermat's Last Theorem that Germain
outlined in her 1819 letter to Gauss, the relationship of these claims to
Sophie Germain's Theorem, and on presenting a coherent and comprehensive
mathematical picture of the many facets of Germain's overall plan of attack on
Fermat's Last Theorem, distilled from the various manuscripts.

We will explain some of Germain's most important mathematical devices in her
approach to Fermat's Last Theorem, provide a sense for the results she
successfully obtained and the ones that are problematic, compare with the
impression of her work left by Legendre's treatise, and in particular discuss
possible overlap between Germain's work and Legendre's. We will also find
connections between Germain's work on Fermat's Last Theorem and that of
mathematicians of the later nineteenth and twentieth centuries. Finally, we
will discuss claims in Germain's manuscripts to have actually fully proven
Fermat's Last Theorem for certain exponents.

Our assessment is based on analyzing all of the following, to which we have
given short suggestive names for reference throughout the paper:

\begin{itemize}
\item \textbf{Manuscript A} (Biblioth\`{e}que Nationale): An undated
manuscript entitled \emph{Remarques sur l'impossibilit\'{e} de satisfaire en
nombres entiers a l'\'{e}quation }$x^{p}+y^{p}=z^{p}$ \cite[pp.~198r--208v]%
{ger9114} (20 sheets numbered in Germain's hand, with 13 carefully labeled
sections). This is a highly polished version of Manuscript D (some, but not
all, of the marginal notes added to Manuscript A have been noted in the
transcription of Manuscript D in \cite{delcentina-germain-flt});

\item \textbf{Errata to Manuscript A }(Biblioth\`{e}que Nationale): Two
undated sheets \cite[pp.~214r, 215v]{ger9114} titled \textquotedblleft
errata\textquotedblright\ or \textquotedblleft erratu\textquotedblright;

\item \textbf{Manuscript B} (Biblioth\`{e}que Nationale): An undated
manuscript entitled \emph{D\'{e}monstration de l'impossibilit\'{e} de
satisfaire en nombres entiers \`{a} l'\'{e}quation }$z^{2(8n\pm3)}%
=y^{2(8n\pm3)}+x^{2(8n\pm3)}$ \cite[pp.~92r--94v]{ger9114} (4 sheets);

\item \textbf{Manuscript C }(Biblioth\`{e}que Nationale): A polished undated
set of three pages \cite[pp.\ 348r--349r]{ger9115} stating and claiming a
proof of Fermat's Last Theorem for all even exponents;

\item \textbf{Letter from Germain to Legendre} (New York Public Library): An
undated 3 page letter \cite{gerleg}\footnote{Although we have found nothing
else in the way of correspondence between Legendre and Germain on Fermat's
Last Theorem, we are fortunate to know of this one critical letter, held in
the Samuel Ward papers of the New York Public Library. These papers include,
according to the collection guide to the papers, \textquotedblleft letters by
famous mathematicians and scientists acquired by Ward with his purchase of the
library of mathematician A. M. Legendre.\textquotedblright\ We thank Louis
Bucciarelli for providing us with this lead.} about Fermat's Last Theorem;

\item \textbf{Manuscript D} (Biblioteca Moreniana): A less polished version of
Manuscript A \cite[cass. 11, ins. 266]{gernfl} \cite[p.\ 234]{delcentinabook}
(25 pages, the 19th blank), transcribed in \cite{delcentina-germain-flt};

\item \textbf{Letter of May 12, 1819 from Germain to Gauss}
(Nieders\"{a}chsische Staats- und Universit\"{a}tsbibliothek G\"{o}ttingen): A
letter of eight numbered sheets \cite{gerauss}, mostly about her work on
Fermat's Last Theorem, transcribed in \cite{delcentina-germain-flt}.
\end{itemize}

Together these appear to be Germain's primary pieces of work on Fermat's Last
Theorem. Nevertheless, our assessment is based on only part of her
approximately 150--200 pages of number theory manuscripts in the
Biblioth\`{e}que, and other researchers may ultimately have more success than
we at deciphering, understanding, and interpreting them. Also, there are
numerous additional Germain papers in the Biblioteca Moreniana that may yield
further insight. Finally, even as our analysis and evaluation answers many
questions, it will also raise numerous new ones, so there is fertile ground
for much more study of her manuscripts by others. In particular, questions of
the chronology of much of her work, and of her interaction with others, still
contain enticing perplexities.

Before beginning our analysis of Germain's manuscripts, we summarize for
comparison Andrea Del Centina's recent work \cite{delcentina-germain-flt}. He
first analyzes an appendix\footnote{Held in the Biblioteca Moreniana.} to an
1804 letter from Germain to Gauss (for which he provides a transcription in
his own appendix). This represents her very early work on Fermat's Last
Theorem, in which she claims (incorrectly) a proof for a certain family of
exponents; this 1804 approach was mathematically unrelated to the coherent
theory that we will see in all her much later manuscripts. Then Del Centina
provides an annotated transcription of the entire 1819 letter to Gauss, which
provides her own not too technical overview for Gauss of her later and more
mature mathematical approach. We focus on just a few translated excerpts from
this 1819 letter, to provide an overview and to introduce key aspects of her
various manuscripts.

Finally Del Centina leads the reader through an analysis of the mathematics in
Manuscript D (almost identical with A), which he also transcribes in its
entirety in an appendix. Although Manuscript A is our largest and most
polished single source, we view it within the context of all the other
manuscripts and letters listed above, since our aim is to present most of
Germain's web of interconnected results in one integrated mathematical
framework, illustrated in Figures \ref{figure-1} and \ref{figure-2} in Section
\ref{precis}. Also, even in the analysis of the single Manuscript A that is
discussed in both Del Centina's paper and ours, we and Del Centina very often
place our emphases on different aspects, and draw somewhat different
conclusions about parts of the manuscript. We will not remark specially on
numerous aspects of Manuscript A that are discussed either only in his paper
or only in ours; the reader should consult both. Our footnotes will largely
comment on differences in the treatment of aspects discussed in both
papers.\footnote{In particular, in section \ref{mistake} we examine a subtle
but critical mistake in Germain's proof of a major result, and her later
attempts to remedy it. In his analysis of the same proof, Del Centina does not
appear to be aware of this mistake or its consequences.} Del Centina does not
mention Germain's Errata to Manuscript A (noted by her in its margin), nor
Manuscripts B or C, or the letter from Germain to Legendre, all of which play
a major role for us.

\subsection{Outline for our presentation of Germain's work}

In Section \ref{gauss} we will examine the interaction and mutual influences
between Germain and Gauss, focusing on Fermat's Last Theorem. In particular we
will display Germain's summary explanation to Gauss in 1819 of her
\textquotedblleft grand plan\textquotedblright\ for proving the impossibility
of the Fermat equation outright, and her description of related successes and
failures. This overview will serve as introduction for reading her main
manuscripts, and to the big picture of her body of work.

The four ensuing Sections \ref{grandplan}, \ref{largesize}, \ref{specialform},
and \ref{evenexponents} contain our detailed analysis of the essential
components of Germain's work. Her mathematical aims included a number of
related results on Fermat's Last Theorem, namely her grand plan, large size of
solutions, $p^{2}$-divisibility of solutions (i.e., Sophie Germain's Theorem,
applicable to Case 1), and special forms of the exponent. These results are
quite intertwined in her manuscripts, largely because the hypotheses that
require verification overlap. We have separated our exposition of these
results in the four sections in a particular way, explained below, partly for
clarity of the big picture, partly to facilitate direct comparison with
Legendre's treatise, which had a different focus but much apparent overlap
with Germain's, and partly to enable easier comparison with the later work of
others. The reader may refer throughout the paper to Figures \ref{figure-1}
and \ref{figure-2} in Section \ref{precis}, which portray the big picture of
the interconnections between Germain's claims (theorems), conditions
(hypotheses), and propositions and algorithms for verifying these conditions.

Section \ref{grandplan} will address Germain's grand plan. We will elucidate
from Manuscripts A and D the detailed methods Germain developed in her grand
plan, the progress she made, and its difficulties. We will compare Germain's
methods with her explanation and claims to Gauss, and with Legendre's work.
The non-consecutivity condition on $p$-th power residues modulo an auxiliary
prime $\theta$, which we saw above in the statement of Sophie Germain's
Theorem, is key also to Germain's grand plan. It has been pursued by later
mathematicians all the way to the present day, and we will compare her
approach to later ones. We will also explore whether Germain at some point
realized that her grand plan could not be carried through, using the published
historical record and a single relevant letter from Germain to Legendre.

Section \ref{largesize} will explore large size of solutions and $p^{2}%
$-divisibility of solutions. In Manuscripts A and D Germain proved and applied
a theorem which we shall call \textquotedblleft Large size of
solutions\textquotedblright, whose intent is to convince that any solutions
which might exist to a Fermat equation would have to be astronomically large,
a claim we will see she mentioned to Gauss in her 1819 letter. Germain's
effort here is challenging to evaluate, since her proof as given in the
primary manuscript is flawed, but she later recognized this and attempted to
compensate. Moreover Legendre published similar results and applications,
which we will contrast with Germain's. We will discover that the theorem on
$p^{2}$-divisibility of solutions that is known in the literature as Sophie
Germain's Theorem is simply minor fallout from her \textquotedblleft Large
size of solutions\textquotedblright\ analysis. And we will compare the methods
she uses to apply her theorem with the methods of later researchers.

Section \ref{specialform} addresses a large family of prime exponents for the
Fermat equation. In Manuscript B, Germain claims proof of Fermat's Last
Theorem for this family of exponents, building on an essentially
self-contained statement of Sophie Germain's Theorem on $p^{2}$-divisibility
of solutions to deal with Case 1 for all exponents first.

Section \ref{evenexponents} considers even exponents. Germain's Manuscript C,
using a very different approach from the others, claims to prove Fermat's Last
Theorem for all even exponents based on the impossibility of another
Diophantine equation.

We end the paper with three final sections: a pr\'{e}cis and connections for
Germain's various thrusts at Fermat's Last Theorem, our reevaluation, and a
conclusion. The reevaluation will take into account Germain's frontal assault
on Fermat's Last Theorem, her analysis leading to claims of astronomical size
for any possible solutions to the Fermat equation, the fact that Sophie
Germain's Theorem is in the end a small piece of something much more
ambitious, our assessment of how independent her work actually was from her
mentor Legendre's, of the methods she invented for verifying various
conditions, and of the paths unknowingly taken in her footsteps by later
researchers. We will conclude that a substantial elevation of Germain's
contribution is in order.

\section{Interactions with Gauss on number theory \label{gauss}}

Number theory held a special fascination for Germain throughout much of her
life. Largely self-taught, due to her exclusion as a woman from higher
education and normal subsequent academic life, she had first studied
Legendre's \emph{Th\'{e}orie des Nombres}, published in 1798, and then
devoured Gauss's \emph{Disquisitiones Arithmeticae} when it appeared in 1801
\cite{libri-notice}. Gauss's work was a complete departure from everything
that came before, and organized number theory as a mathematical subject
\cite{goldstein} \cite{neumann}, with its own body of methods, techniques, and
objects, including the theory of congruences and the roots of the cyclotomic equation.

\subsection{Early correspondence}

Germain's exchange of letters with Gauss, initiated under the male pseudonym
LeBlanc, lasted from 1804 to 1808, and gave tremendous impetus to her work. In
her first letter \cite{boncompagni}\footnote{Relevant excerpts can be found in
Chapter 3 of \cite{bucc}; see also \cite{stupuy}.} she sent Gauss some initial
work on Fermat's Last Theorem stemming from inspiration she had received from
his \emph{Disquisitiones}.

Gauss was greatly impressed by Germain's work, and was even stimulated thereby
in some of his own, as evidenced by his remarks in a number of letters to his
colleague Wilhelm Olbers. On September 3, 1805 Gauss wrote \cite[p.~268]%
{olbers}: \textquotedblleft Through various circumstances --- partly through
several letters from LeBlanc in Paris, who has studied my \emph{Disq. Arith.}
with a true passion, has completely mastered them, and has sent me occasional
very respectable communications about them, [$\ldots$] I have been tempted
into resuming my beloved arithmetic investigations.\textquotedblright%
\footnote{Throughout the paper, any English translations are our own, unless
cited otherwise.}\ After LeBlanc's true identity was revealed to him, he wrote
again to Olbers, on March 24, 1807 \cite[p.~331]{olbers}: \textquotedblleft
Recently my \emph{Disq. Arith.} caused me a great surprise. Have I not written
to you several times already about a correspondent LeBlanc from Paris, who has
given me evidence that he has mastered completely all investigations in this
work? This LeBlanc has recently revealed himself to me more closely. That
LeBlanc is only a fictitious name of a young lady Sophie Germain surely amazes
you as much as it does me.\textquotedblright\ 

Gauss's letter to Olbers of July 21 of the same year shows that Germain had
become a valued member of his circle of correspondents \cite[pp.~376--377]%
{olbers}: \textquotedblleft Upon my return I\ have found here several letters
from Paris, by Bouvard, Lagrange, and Sophie Germain. [$\ldots$] Lagrange
still shows much interest in astronomy and higher arithmetic; the two sample
theorems (for which prime numbers\footnote{as modulus.} is [the number] two a
cubic or biquadratic residue), which I also told you about some time ago, he
considers `that which is most beautiful and difficult to prove.' But Sophie
Germain has sent me the proofs for them; I have not yet been able to look
through them, but I\ believe they are good; at least she has approached the
matter from the right point of view, only they are a little more long-winded
than will be necessary.\textquotedblright\ 

The two theorems on power residues were part of a letter Gauss wrote to
Germain on April 30, 1807 \cite[vol.~10,~pp.~70--74]{gausswerke}. Together
with these theorems he also included, again without proof, another result now
known as Gauss's Lemma, from which he says one can derive special cases of the
Quadratic Reciprocity Theorem, the first deep result discovered and proven
about prime numbers.\footnote{Gauss was the first to prove quadratic
reciprocity, despite major efforts by both its discoverer Euler and by
Legendre.} In a May 12, 1807 letter to Olbers, Gauss says \textquotedblleft
Recently I replied to a letter of hers and shared some Arithmetic with her,
and this led me to undertake an inquiry again; only two days later I\ made a
very pleasant discovery. It is a new, very neat, and short proof of the
fundamental theorem of art. 131.\textquotedblright\ \cite[pp.~360]{olbers} The
proof Gauss is referring to, based on the above lemma in his letter to
Germain, is now commonly called his \textquotedblleft third\textquotedblright%
\ proof of the Quadratic Reciprocity Theorem, and was published in 1808
\cite{gauss1808}, where he says he has finally found \textquotedblleft the
simplest and most natural way to its proof\textquotedblright\ (see also
\cite{lp1,lp2}).

We shall see in Germain's manuscripts that the influence of Gauss's
\emph{Disquisitiones} on her work was all-encompassing; her manuscripts and
letters use Gauss's congruence notion and point of view throughout, in
contrast to her Paris mentor Legendre's style of equalities \textquotedblleft
omitting multiples\textquotedblright\ of the modulus. Her work benefits from
the ease of writing and thinking in terms of arithmetic modulo a prime enabled
by the \emph{Disquisitiones} \cite{goldstein} \cite{neumann,wussing}. Germain
also seems to have been one of the very first to adopt and internalize in her
own research the ideas of the \emph{Disquisitiones}. But her work, largely
unpublished, may have had little influence on the next generation.

\subsection{Letter of 1819 about Fermat's Last Theorem}

On the twelfth of May, 1819, Sophie Germain penned a letter from her Parisian
home to Gauss in G\"{o}ttingen \cite{gerauss}. Most of this lengthy letter
describes her work on substantiating Fermat's Last Theorem.

The letter provides a window into the context of their interaction on number
theory from a vantage point fifteen years after their initial correspondence.
It will show us how she viewed her overall work on Fermat's Last Theorem at
that time, placing it in the bigger picture of her mathematical research, and
specifically within her interaction with and influence from Gauss. And the
letter will give enough detail on her actual progress on proving Fermat's Last
Theorem to prepare us for studying her manuscripts, and to allow us to begin
comparison with the published historical record, namely the attribution by
Legendre in 1823 of Sophie Germain's Theorem.%

\begin{figure}
[pt]
\begin{center}
\includegraphics[
natheight=7.675200in,
natwidth=7.203900in,
height=5.3748in,
width=5.047in
]%
{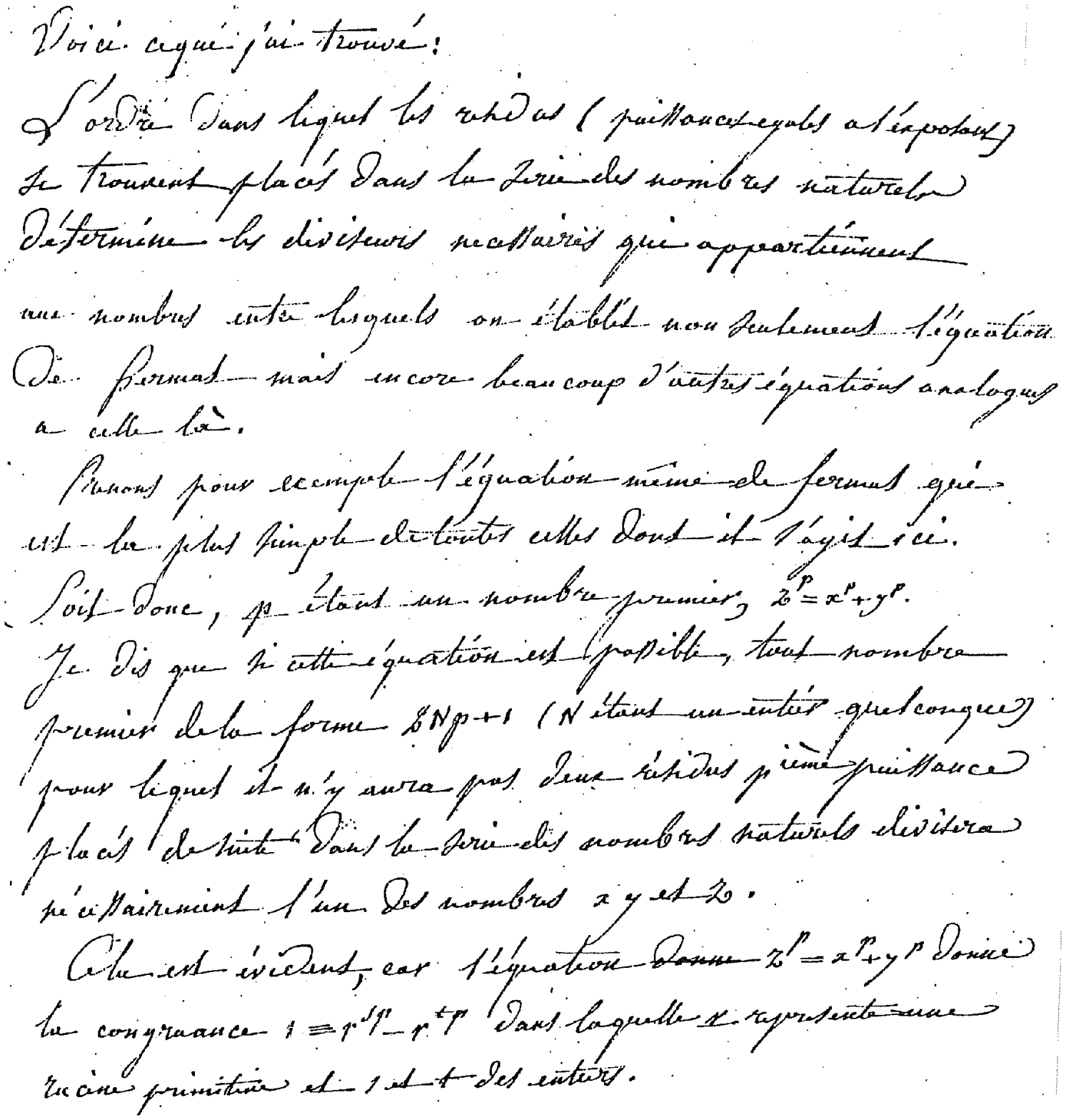}%
\caption{\textquotedblleft Voici ce que ja'i trouv\'{e}:\textquotedblright%
\ From Germain's letter to Gauss, 1819}%
\label{gauss-l-image}%
\end{center}
\end{figure}

Germain's letter was written after an eleven year hiatus in their
correspondence. Gauss had implied in his last letter to Germain in 1808 that
he might not continue to correspond due to his new duties as astronomer, but
the visit of a friend of Gauss's to Paris in 1819 provided Germain the
encouragement to attempt to renew the exchange \cite[p. 86, 137]{bucc}. She
had a lot to say. Germain describes first the broad scope of many years of
work, to be followed by details on her program for proving Fermat's Last Theorem:

\begin{quotation}
[...] Although I have worked for some time on the theory of vibrating surfaces
[...], I have never ceased thinking about the theory of numbers. I will give
you a sense of my absorption with this area of research by admitting to you
that even without any hope of success, I still prefer it to other work which
might interest me while I think about it, and which is sure to yield results.

Long before our Academy proposed a prize for a proof of the impossibility of
the Fermat equation, this type of challenge, which was brought to modern
theories by a geometer who was deprived of the resources we possess today,
tormented me often. I glimpsed vaguely a connection between the theory of
residues and the famous equation; I believe I spoke to you of this idea a long
time ago, because it struck me as soon as I read your
book.\footnote{\textquotedblleft Quoique j'ai travaill\'{e} pendant quelque
tems a la th\'{e}orie des surfaces vibrantes [$\dots$], je n'ai jamais
cess\'{e} de penser a la th\'{e}orie des nombres. Je vous donnerai une
id\'{e}e de ma pr\'{e}occupation pour ce genre de recherches en vous avouant
que m\^{e}me sans aucune esperance de succ\`{e}s je la prefere a un travail
qui me donnerait necessairement un resultat et qui pourtant m'interresse
$\dots$ quand j'y pense.
\par
\textquotedblleft Longtems avant que notre academie ait propos\'{e} pour sujet
de prix la d\'{e}monstration de l'impossibilit\'{e} de l'\'{e}quation de
Fermat cet espece de d\'{e}fi---port\'{e} aux th\'{e}ories modernes par un
g\'{e}ometre --- qui fut priv\'{e} des resources que nous possedons
aujourd'hui me tourmentait souvent. J'entrevoyais \emph{vaguement} une liaison
entre la th\'{e}orie des residus et la fameuse \'{e}quation, je crois m\^{e}me
vous avoir parl\'{e} anciennement de cette id\'{e}e car elle m'a frapp\'{e}
aussit\^{o}t que j'ai connu votre livre.\textquotedblright\ (Letter to Gauss,
p. 2)}
\end{quotation}

Germain continues the letter by explaining to Gauss her major effort to prove
Fermat's Last Theorem (Figure \ref{gauss-l-image}), including the overall
plan, a summary of results, and claiming to have proved the astronomically
large size of any possible solutions. She introduces her work to him with the
words \textquotedblleft Voici ce que ja'i trouv\'{e}:\textquotedblright%
\ (\textquotedblleft Here is what I have found:\textquotedblright).

\begin{quotation}
Here is what I have found: [...]

The order in which the residues (powers equal to the exponent\footnote{i.e.,
power residues where the power is equal to the exponent in the Fermat
equation.}) are distributed in the sequence of natural numbers determines the
necessary divisors which belong to the numbers among which one establishes not
only the equation of Fermat, but also many other analogous equations.

Let us take for example the very equation of Fermat, which is the simplest of
those we consider here. Therefore we have $z^{p}=x^{p}+y^{p}$, $p$ a prime
number. I claim that if this equation is possible, then every prime number of
the form $2Np+1$ ($N$ being any integer), for which there are no two
consecutive $p$-th power residues in the sequence of natural
numbers,\footnote{\label{thetaform}Germain is considering congruence modulo an
auxiliary prime $\theta=2Np+1$ that has no consecutive nonzero $p$-th power
residues. While the specified form of $\theta$ is not necessary to her
subsequent argument, she knows that only prime moduli of the form
$\theta=2Np+1$ can possibly have no consecutive nonzero $p$-th power residues,
and implicitly that Gauss will know this too. (This is easy to confirm using
Fermat's ``Little'' Theorem; see, for instance, \cite[p.\ 124]%
{ribenboim-amateurs}.) Thus she restricts without mention to considering only
those of this form.} necessarily divides one of the numbers $x$, $y$, and $z$.

This is clear, since the equation $z^{p}=x^{p}+y^{p}$ yields the congruence
$1\equiv r^{sp}-r^{tp}$ in which $r$ represents a primitive root and $s$ and
$t$ are integers.\footnote{Here Germain is utilizing two facts about the
residues modulo the prime $\theta$. One is that when the modulus is prime, one
can actually \textquotedblleft divide\textquotedblright\ in modular arithmetic
by any number with nonzero residue. So if none of $x,y,z$ were divisible by
$\theta$, then modular division of the Fermat equation by $x^{p}$ or $y^{p}$
would clearly produce two nonzero consecutive $p$-th power residues. She is
also using the fact that for a prime modulus, there is always a number, called
a primitive root for this modulus, such that any number with nonzero residue
is congruent to a power of the primitive root. She uses this representation in
terms of a primitive root later on in her work.}\ [...]

It follows that if there were infinitely many such numbers, the equation would
be impossible.

I have never been able to arrive at the infinity, although I have pushed back
the limits quite far by a method of trials too long to describe here. I still
dare not assert that for each value of $p$ there is no limit beyond which all
numbers of the form $2Np+1$ have two consecutive $p$-th power residues in the
sequence of natural numbers. This is the case which concerns the equation of Fermat.

You can easily imagine, Monsieur, that I have been able to succeed at proving
that this equation is not possible except with numbers whose size frightens
the imagination; because it is also subject to many other conditions which I
do not have the time to list because of the details necessary for establishing
its success. But all that is still not enough; it takes the infinite and not
merely the very large.\footnote{\textquotedblleft Voici ce que j'ai trouv\'{e}
:
\par
\textquotedblleft L'ordre dans lequel les residus (puissances egales a
l'exposant) se trouvent plac\'{e}s dans la serie des nombres naturels
d\'{e}termine les diviseurs necessaires qui appartiennent aux nombres entre
lequels on \'{e}tablit non seulement l'\'{e}quation de Fermat mais encore
beaucoup d'autres \'{e}quations analogues a celle l\`{a}.
\par
\textquotedblleft Prenons pour exemple l'\'{e}quation m\^{e}me de Fermat qui
est la plus simple de toutes celles dont il s'agit ici. Soit donc, $p$
\'{e}tant un nombre premier, $z^{p}=x^{p}+y^{p}$. Je dis que si cette
\'{e}quation est possible, tout nombre premier de la forme $2Np+1$ ($N$
\'{e}tant un entier quelconque) pour lequel il n'y aura pas deux r\'{e}sidus
$p^{\text{i\`{e}me}}$ puissance plac\'{e}s de suite dans la serie des nombres
naturels divisera n\'{e}cessairement l'un des nombres $x$ $y$ et $z$.
\par
\textquotedblleft Cela est \'{e}vident, car l'\'{e}quation $z^{p}=x^{p}+y^{p}$
donne la congruence $1\equiv r^{sp}-r^{tp}$ dans laquelle $r$ represente une
racine primitive et $s$ et $t$ des entiers.
\par
\textquotedblleft$\ldots$ Il suit del\`{a} que s'il y avoit un nombre infini
de tels nombres l'\'{e}quation serait impossible.
\par
\textquotedblleft Je n'ai jamais p\^{u} arriver a l'infini quoique j'ai
recul\'{e} bien loin les limites par une methode de tatonnement trop longue
pour qu'il me soit possible de l'exposer ici. Je n'oserais m\^{e}me pas
affirmer que pour chaque valeur de $p$ il n'existe pas une limite audela
delaquelle tous les nombres de la forme $2Np+1$ auraient deux r\'{e}sidus
$p^{\text{i\`{e}mes}}$ plac\'{e}s de suite dans la serie des nombres naturels.
C'est le cas qui interesse l'\'{e}quation de Fermat.
\par
\textquotedblleft Vous concevrez aisement, Monsieur, que j'ai d\^{u} parvenir
a prouver que cette \'{e}quation ne serait possible qu'en nombres dont la
grandeur effraye l'imagination ; Car elle est encore assujettie a bien
d'autres conditions que je n'ai pas le tems d'\'{e}num\'{e}rer a cause des
details necessaire pour en \'{e}tablir la r\'{e}ussite. Mais tout cela n'est
encore rien, il faut l'infini et non pas le tr\`{e}s grand.\textquotedblright%
\ (Letter to Gauss, pp. 2--4)}
\end{quotation}

Several things are remarkable here. Most surprisingly, Germain does not
mention to Gauss anything even hinting at the only result she is actually
known for in the literature, what we call Sophie Germain's Theorem. Why not?
Where is it? Instead, Germain explains a plan, simple in its conception, for
proving Fermat's Last Theorem outright. It requires that, for a given prime
exponent $p$, one establish infinitely many auxiliary primes each satisfying a
non-consecutivity condition on its nonzero $p$-th power residues (note that
this condition is the very same as one of the two hypotheses required in
Sophie Germain's Theorem for proving Case 1, but there one only requires a
single auxiliary prime, not infinitely many). And she explains to Gauss that
since each such auxiliary prime will have to divide one of $x$, $y$, $z$, the
existence of infinitely many of them will make the Fermat equation impossible.
She writes that she has worked long and hard at this plan by developing a
method for verifying the condition, made great progress, but has not been able
to bring it fully to fruition (even for a single $p$) by verifying the
condition for infinitely many auxiliary primes. She also writes that she has
proven that any solutions to a Fermat equation would have to \textquotedblleft
frighten the imagination\textquotedblright\ with their size. And she gives a
few details of her particular methods of attack. The next two sections will
examine the details of these claims in Germain's manuscripts.

\section{The grand plan\label{grandplan}}%

\begin{figure}
[pt]
\begin{center}
\includegraphics[
natheight=3.578600in,
natwidth=6.578600in,
height=2.7579in,
width=5.047in
]%
{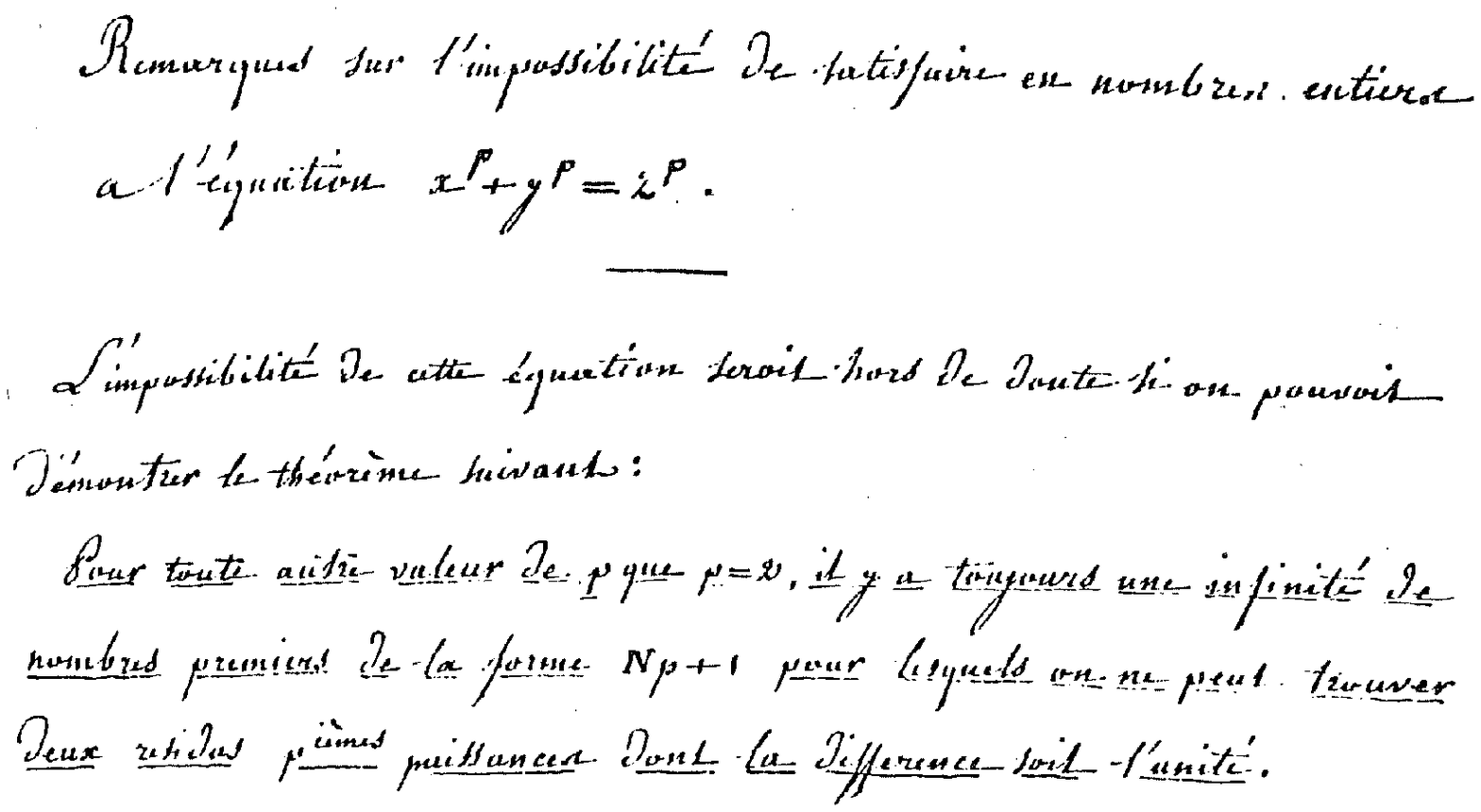}%
\caption{Beginning of Manuscript A}%
\label{man-a-198top-image}%
\end{center}
\end{figure}

Our aim in this section is to study Germain's plan for proving Fermat's Last
Theorem, as outlined to Gauss, to show its thoroughness and sophistication,
and to consider its promise for success.

As we saw Germain explain to Gauss, one can prove Fermat's Last Theorem for
exponent $p$ by producing an infinite sequence of qualifying auxiliary primes.
Manuscript A (Figure \ref{man-a-198top-image}) contains, among other things,
the full details of her efforts to carry this plan through, occupying more
than 16 pages of very polished writing. We analyze these details in this
section, ending with a comparison between Manuscripts A and D.

\subsection{Germain's plan for proving Fermat's Last
Theorem\label{germain's plan}}

We have seen that Germain's plan for proving Fermat's Last Theorem for
exponent $p$ hinged on developing methods to validate the following qualifying
condition for infinitely many auxiliary primes of the form $\theta=2Np+1$:

\begin{NC}
[\textbf{N}on-\textbf{C}onsecutivity]There do \textbf{n}ot exist two nonzero
\textbf{c}onsecutive $p^{\text{th}}$ power residues, modulo $\theta$.
\end{NC}

%

\begin{figure}
[pt]
\begin{center}
\includegraphics[
natheight=5.090300in,
natwidth=7.087100in,
height=3.6331in,
width=5.047in
]%
{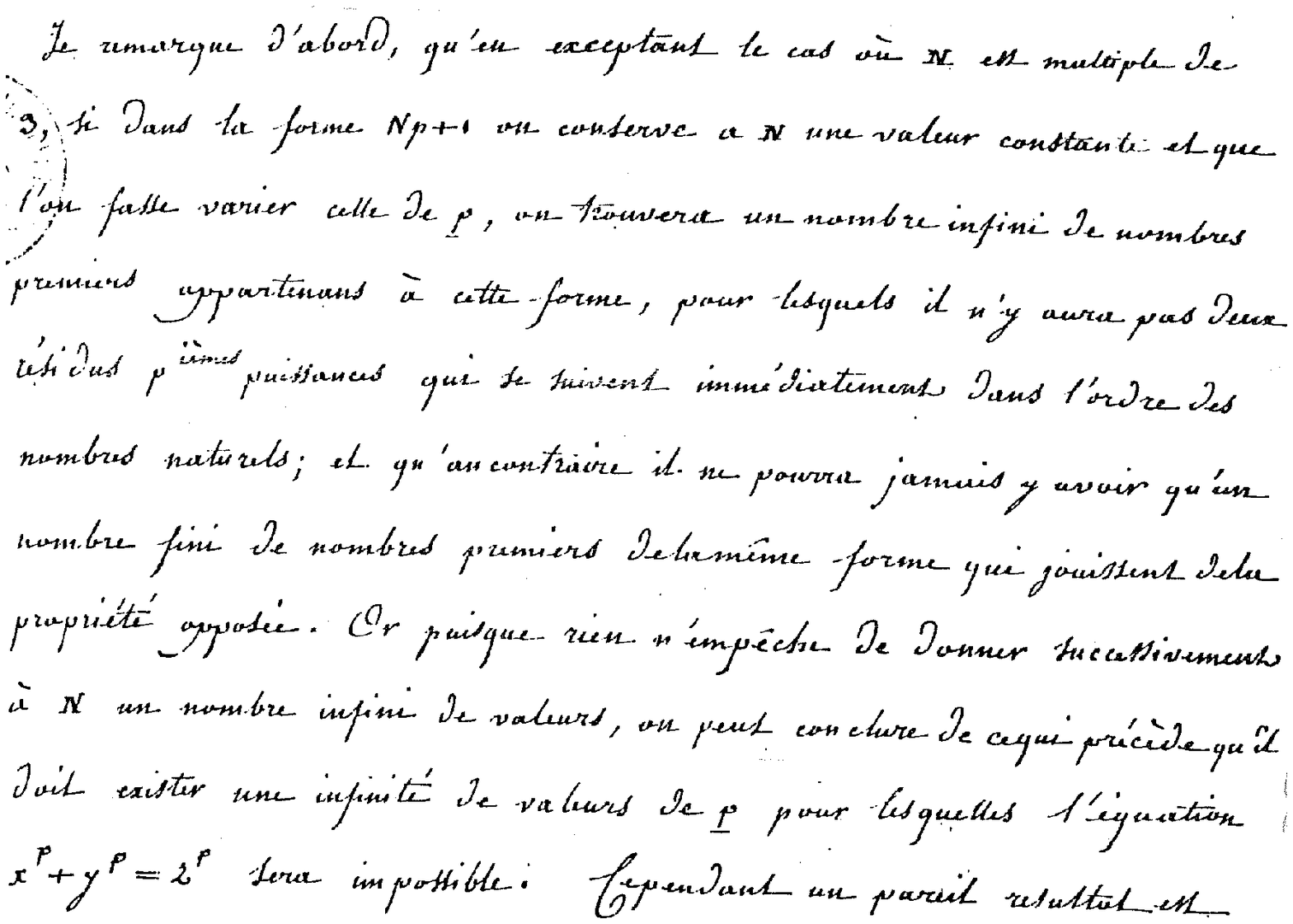}%
\caption{From the introduction of Manuscript A}%
\label{man-a-198bot-image}%
\end{center}
\end{figure}

Early on in Manuscript A (Figure \ref{man-a-198bot-image}), Germain claims
that for each fixed $N$ (except when $N$ is a multiple of $3,$ for which she
shows that Condition N-C always fails\footnote{See \cite[p.\ 127]%
{ribenboim-amateurs}.}), there will be only finitely many exceptional numbers
$p$ for which the auxiliary $\theta=2Np+1$ fails to satisfy Condition N-C
(recall from footnote \ref{thetaform} that only primes of the form
$\theta=2Np+1$ can possibly satisfy the N-C condition). Much of Germain's
manuscript is devoted to supporting this claim; while she was not able to
bring this to fruition, Germain's insight was vindicated much later when
proven true by E. Wendt in 1894 \cite[p. 756]{dickson} \cite[p. 124ff]%
{ribenboim-amateurs} \cite{wendt}.\footnote{Germain's claim would follow
immediately from Wendt's recasting of the condition in terms of a circulant
determinant depending on $N$: Condition N-C fails to hold for $\theta$ only if
$p$ divides the determinant, which is nonzero for all $N$ not divisible by
$3$. There is no indication that Wendt was aware of Germain's work.}

Note that a priori there is a difference in impact between analyzing Condition
N-C for fixed $N$ versus for fixed $p$. To prove Fermat's Last Theorem for
fixed $p$, one needs to verify N-C for infinitely many $N$, whereas Germain's
approach is to fix $N$ and aim to verify N-C for all but finitely many $p$.
Germain was acutely aware of this distinction. After we see exactly what she
was able to accomplish for fixed $N$, we will see what she had to say about
converting this knowledge into proving Fermat's Last Theorem for particular
values of $p$.

Before delving into Germain's reasoning for general $N$, let us consider just
the case $N=1$, i.e., when $\theta=2p+1$ is also prime, today called a
\textquotedblleft Germain prime\textquotedblright. We consider $N=1$ partly
because it is illustrative and not hard, and partly to relate it to the
historical record. Germain knew well that there are always precisely $2N$
nonzero $p$-th power residues modulo an auxiliary prime of the form
$\theta=2Np+1$. Thus in this case, the numbers $1$ and $2p=\theta-1\equiv-1$
are clearly the only nonzero $p$-th power residues, so Condition N-C
automatically holds. Of course for $N>1$, with more $p$-th power residues,
their distribution becomes more difficult to analyze. Regarding the historical
record, we remark that the other condition of Sophie Germain's Theorem\ for
Case 1, namely that $p$ itself not be a $p$-th power modulo $\theta$, is also
obviously satisfied in this case. So Sophie Germain's Theorem automatically
proves Case 1 whenever $2p+1$ is prime. This may shed light on why, as
mentioned earlier, some writers have incorrectly thought that Sophie Germain's
Theorem\ deals only with Germain primes as auxiliaries.

\subsubsection{Establishing Condition N-C for each $N$, including an induction
on $N$}

In order to establish Condition N-C for various $N$ and $p$, Germain engages
in extensive analysis over many pages of the general consequences of nonzero
consecutive $p$-th power residues modulo a prime $\theta=2Np+1$ ($N$ never a
multiple of $3$).

Her analysis actually encompasses all natural numbers for $p$, not just
primes. This is important in relation to the form of $\theta$, since she
intends to carry out a mathematical induction on $N$, and eventually explains
in detail her ideas about how the induction should go. She employs throughout
the notion and notation of congruences introduced by Gauss, and utilizes to
great effect a keen understanding that the $2Np$ multiplicative units mod
$\theta$ are cyclic, generated by a primitive $2Np$-th root of unity, enabling
her to engage in detailed analyses of the relative placement of the nonzero
$p$-th powers (i.e., the $2N$-th roots of $1$) amongst the residues. She is
acutely aware (expressed by us in modern terms) that subgroups of the group of
units are also cyclic, and of their orders and interrelationships, and uses
this in a detailed way. Throughout her analyses she deduces that in many
instances the existence of nonzero consecutive $p$-th power residues would
ultimately force $2$ to be a $p$-th power mod $\theta$, and she therefore
repeatedly concludes that Condition N-C holds under the following hypothesis:

\begin{2Np}
[$2$ is \textbf{N}ot a $p$-th power]The number $2$ is \textbf{n}ot a $p$-th
power residue, modulo $\theta$.
\end{2Np}

Notice that this hypothesis is always a necessary condition for Condition N-C
to hold, since if $2$ is a $p$-th power, then obviously $1$ and $2$ are
nonzero consecutive $p$-th powers; so making this assumption is no
restriction, and Germain is simply exploring whether $2$-N-$p$ is also
sufficient to ensure N-C.

Always assuming this hypothesis, whose verification we shall discuss in
Section \ref{sec2np}, and also the always necessary condition mentioned above
(Section \ref{germain's plan}) that $N$ is not a multiple of $3$, Germain's
analysis initially shows that if there exist two nonzero consecutive $p$-th
power residues, then by inverting them, or subtracting them from $-1$, or
iterating combinations of these transformations, she can obtain more pairs of
nonzero consecutive $p$-th power residues.\footnote{In fact these
transformations are permuting the pairs of consecutive residues according to
an underlying group with six elements, which we shall discuss later. Germain
even notes, when explaining the situation in her letter to Gauss
\cite{gerauss}, that from any one of the six pairs, her transformations will
reproduce the five others.}

Germain proves that, under her constant assumption that $2$ is not a $p$-th
power residue modulo $\theta$, this transformation process will produce at
least $6$ completely disjoint such pairs, i.e., involving at least $12$ actual
$p$-th power residues.\footnote{Del Centina \cite[p. 367ff]%
{delcentina-germain-flt} provides details of how Germain proves this.}
Therefore since there are precisely $2N$ nonzero $p$-th power residues modulo
$\theta$, she instantly proves Condition N-C for all auxiliary primes $\theta$
with $N=1,2,4,5$ as long as $p$ satisfies Condition $2$-N-$p$. Germain
continues with more detailed analysis of these permuted pairs of consecutive
$p$-th power residues (still assuming Condition $2$-N-$p$) to verify Condition
N-C for $N=7$ (excluding $p=2$) and $N=8$ (here she begins to use inductive
information for earlier values of $N$).\footnote{Del Centina \cite[p.
369ff]{delcentina-germain-flt} provides details for $N=7,$ $8$.}

At this point Germain explains her general plan to continue the method of
analysis to higher $N$, and how she would use induction on $N$ for all $p$
simultaneously. In a nutshell, she argues that the existence of nonzero
consecutive $p$-th power residues would have to result in a pair of nonzero
consecutive $p$-th powers, $x,$ $x+1$, for which $x$ is (congruent to) an odd
power (necessarily less than $2N$) of $x+1$. She claims that one must then
analyze cases of the binomial expansion of this power of $x+1$, depending on
the value of $N$, to arrive at the desired contradiction, and she carries out
a complete detailed calculation for $N=10$ (excluding $p=2,3$) as a specific
\textquotedblleft example\textquotedblright\footnote{(Manuscript A, p.
13)}\ of how she says the induction will work in general.\footnote{Del Centina
\cite[p. 369ff]{delcentina-germain-flt} also has commentary on this.}

It is difficult to understand fully this part of the manuscript. Germain's
claims may in fact hold, but we cannot verify them completely from what she
says. Germain's mathematical explanations often omit many details, leaving
much for the reader to fill in, and in this case, there is simply not enough
detail to make a full judgement. Specifically, we have difficulty with an
aspect of her argument for $N=7$, with her explanation of exactly how her
mathematical induction will proceed, and with an aspect of her explanation of
how in general a pair $x,$ $x+1$ with the property claimed above is ensured.
Finally, Germain's example\ calculation for $N=10$ is much more ad hoc than
one would like as an illustration of how things would go in a mathematical
induction on $N$. It seems clear that as this part of the manuscript ends, she
is presenting only a sketch of how things could go, indicated by the fact that
she explicitly states that her approach to induction is via the example of
$N=10$, which is not presented in a way that is obviously generalizable.
Nonetheless, her instincts here were correct, as proven by Wendt.

\subsubsection{The interplay between $N$ and $p$}

Recall from above that proving Condition N-C for all $N$, each with finitely
many excepted $p$, does not immediately solve the Fermat problem.

What is actually needed, for each fixed prime $p$, is that N-C holds for
infinitely many $N$, not the other way around. For instance, perhaps $p=3$
must be excluded from the validation of Condition N-C for all sufficiently
large $N$, in which case Germain's method would not prove Fermat's Last
Theorem for $p=3$. Germain makes it clear early in the manuscript that she
recognizes this issue, that her results do not completely resolve it, and that
she has not proved Fermat's claim for a single predetermined exponent. But she
also states that she strongly believes that the needed requirements do in fact
hold, and that her results for $N\leq10$ strongly support this. Indeed, note
that so far the only odd prime excluded in any verification was $p=3$ for
$N=10$ (recall, though, that we have not yet examined Condition $2$-N-$p$,
which must also hold in all her arguments, and which will also exclude certain
combinations of $N$ and $p$ when it fails).

Germain's final comment on this issue states first that as one proceeds to
ever higher values of $N$, there is always no more than a \textquotedblleft
very small number\textquotedblright\footnote{(Manuscript A, p. 15)}\ of values
of $p$ for which Condition N-C fails. If indeed this, the very crux of the
whole approach, were the case, in particular if the number of such excluded
$p$ were bounded uniformly, say by $K$, for all $N$, which is what she in
effect claims, then a little reflection reveals that indeed her method would
have proven Fermat's Last Theorem for all but $K$ values of $p$, although one
would not necessarily know for which values. She herself then states that this
would prove the theorem for infinitely many $p$, even though not for a single
predetermined value of $p$. It is in this sense that Germain believed her
method could prove infinitely many instances of Fermat's Last Theorem.

\subsubsection{Verifying Condition $2$-N-$p$\label{sec2np}}

We conclude our exposition of Germain's grand plan in Manuscript A with her
subsequent analysis of Condition $2$-N-$p$, which was required for all her
arguments above.

She points out that for $2$ to be a $p$-th power mod $\theta=2Np+1$ means that
$2^{2N}\equiv1$ $\left(  \operatorname{mod}\text{ }\theta\right)  $ (since the
multiplicative structure is cyclic). Clearly for fixed $N$ this can only occur
for finitely many $p$, and she easily determines these exceptional cases
through $N=10$, simply by calculating and factoring each $2^{2N}-1$ by hand,
and observing whether any of the prime factors are of the form $2Np+1$ for any
natural number $p$. To illustrate, for $N=7$ she writes that
\[
2^{14}-1=3\cdot43\cdot127=3\cdot(14\cdot3+1)\cdot\left(  14\cdot9+1\right)  ,
\]
so that $p=3,$ $9$ are the only values for which Condition $2$-N-$p$ fails for
this $N$.

Germain then presents a summary table of all her results verifying Condition
N-C for auxiliary primes $\theta$ using relevant values of $N\leq10$ and
primes $2<p<100$, and says that it can easily be extended
further.\footnote{The table is slightly flawed in that she includes
$\theta=43=14\cdot3+1$ for $N=7$ despite the excluding calculation we just
illustrated, which Germain herself had just written out; it thus seems that
the manuscript may have simple errors, suggesting it may sadly never have
received good criticism from another mathematician.} The results in the table
are impressive. Aside from the case of $\theta=43=14\cdot3+1$ just
illustrated, the only other auxiliary primes in the range of her table which
must be omitted are $\theta=31=10\cdot3+1$, which she determines fails
Condition $2$-N-$p$, and $\theta=61=20\cdot3+1$, which was an exception in her
N-C analysis for $N=10$. In fact each $N$ in her table ends up having at least
five primes $p$ with $2<p<100$ for which $\theta=2Np+1$ is also prime and
satisfies the N-C condition.

While the number of $p$ requiring exclusion for Condition $2$-N-$p$ may appear
\textquotedblleft small\textquotedblright\ for each $N$, there seems no
obvious reason why it should necessarily be uniformly bounded for all $N$;
Germain does not discuss this issue specifically for Condition $2$-N-$p$. As
indicated above, without such a bound it is not clear that this method could
actually prove any instances of Fermat's theorem.

\subsubsection{Results of the grand plan}

As we have seen above, Germain had a sophisticated and highly developed plan
for proving Fermat's Last Theorem for infinitely many exponents.

It relied heavily on facility with the multiplicative structure in a cyclic
prime field and a set (group) of transformations of consecutive $p$-th powers.
She carried out her program on an impressive range of values for the necessary
auxiliary primes, believed that the evidence indicated one could push it
further using mathematical induction by her methods, and she was optimistic
that by doing so it would prove Fermat's Last Theorem for infinitely many
prime exponents. In hindsight we know that, promising as it may have seemed at
the time, the program can never be carried to completion, as we shall see next.

\subsection{Failure of the grand plan}

Did Germain ever know that her grand plan cannot succeed? To answer this
question we examine the published record, Germain's correspondence with Gauss,
and a letter she wrote to Legendre.

Published indication that Germain's method cannot succeed in proving Fermat's
Last Theorem, although not mentioning her by name, came in work of Guglielmo
(Guillaume) Libri, a rising mathematical star in the 1820s. We now describe
Libri's work in this regard.

\subsubsection{Libri's claims that such a plan cannot work}

It is a bit hard to track and compare the content of Libri's relevant works
and their dates, partly because Libri presented or published several different
works all with the same title, but some of these were also multiply published.
Our interest is in the content of just two different works. In 1829 Libri
published a set of his own memoirs \cite{libri-1829}. One of these is titled
\emph{M\'{e}moire sur la th\'{e}orie des nombres}, republished later word for
word as three papers in Crelle's Journal \cite{libri}. The memoir published in
1829 ends by applying Libri's study of the number of solutions of various
congruence equations to the situation of Fermat's Last Theorem. Among other
things, Libri shows that for exponents $3$ and $4$, there can be at most
finitely many auxiliary primes satisfying the N-C condition. And he claims
that his methods will clearly show the same for all higher exponents. Libri
explicitly notes that his result proves that the attempts of others to prove
Fermat's Last Theorem by finding infinitely many such auxiliaries are in vain.

Libri also writes in his 1829 memoir that all the results he obtains were
already presented in two earlier memoirs of 1823 and 1825 to the Academy of
Sciences in Paris. Libri's 1825 presentation to the Academy was also
published, in 1833/1838 \cite{libri-1825}, confusingly with the same title as
the 1829 memoir. This presumably earlier document\footnote{One can wonder when
the document first published in 1833, but based on Libri's 1825 Academy
presentation, was really written or finalized. Remarks he makes in it suggest,
though, that it was essentially his 1825 presentation.} is quite similar to
the publication of 1829, in that it develops methods for determining the
number of solutions to quite general congruence equations, including that of
the N-C condition, but it does not explicitly work out the details for the N-C
condition applying to Fermat's Last Theorem, as did the 1829 memoir.

Thus it seems that close followers of the Academy should have been aware by
1825 that Libri's work would doom the auxiliary prime approach to Fermat's
Last Theorem, but it is hard to pin down exact dates.\footnote{For
completeness, we mention that Libri also published a memoir on number theory
in 1820, his very first publication, with the title \emph{Memoria sopra la
teoria dei numeri} \cite{libri-1820}, but it was much shorter and does not
contain the same type of study or results on the number of solutions to
congruence equations.} Much later, P. Pepin \cite[pp.~318--319]{pepin-1876}
\cite{pepin-1880} and A.-E. Pellet \cite[p.~93]{pellet} (see \cite[p. 750,
753]{dickson} \cite[pp.~292--293]{ribenboim-amateurs}) confirmed all of
Libri's claims, and L. E. Dickson \cite{dickson1909,dickson1909-2} gave
specific bounds.

\subsubsection{What Germain knew and when: Gauss, Legendre, and Libri}

Did Germain ever know from Libri or otherwise that her grand plan to prove
Fermat's Last Theorem could not work, and if so, when?

We know that in 1819 she was enthusiastic in her letter to Gauss about her
method for proving Fermat's Last Theorem, based on extensive work exemplified
by Manuscript A.\footnote{\label{poinsot}Near the end she even expresses to
Gauss how a brand new work by L. Poinsot \cite{poinsot} will help her further
her efforts to confirm the N-C condition by giving a new way of working with
the $p$-th powers mod $\theta=2Np+1$. She interprets them as the solutions of
the binomial equation of degree $2N$, i.e., of $x^{2N}-1=0$. Poinsot's memoir
takes the point of view that the mod $\theta$ solutions of this equation can
be obtained by first considering the equation over the complex numbers, where
much was already known about the complex $2N$-th roots of unity, and then
considering these roots as mod $p$ integers by replacing the complex number
$\sqrt{-1}$ by an integer whose square yields $-1$ mod $p$. Del Centina
\cite[p. 361]{delcentina-germain-flt} also discusses this connection.} In the
letter Germain details several specific examples of her results on the N-C
condition that match perfectly with Manuscript A, and which she explicitly
explains have been extracted from an already much older note
(\textquotedblleft d'une note dej\'{a} ancienne\textquotedblright%
\footnote{(Letter to Gauss, p. 5)}) that she has not had the time to recheck.
In fact everything in the extensive letter to Gauss matches the details of
Manuscript A. This suggests that Manuscript A is likely the older note in
question, and considerably predates her 1819 letter to Gauss. Thus 1819 is our
lower bound for the answer to our question.

We also know that by 1823 Legendre had written his memoir crediting Germain
with her theorem, but without even mentioning the method of finding infinitely
many auxiliary primes that Germain had pioneered to try to prove Fermat's Last
Theorem.\footnote{Del Centina \cite[p. 362]{delcentina-germain-flt} suggests
that a letter from Legendre to Germain in late 1819, published in
\cite{stupuy}, shows that he believed at that time that Germain's work on
Fermat's Last Theorem could not succeed. However, we are not certain that this
letter is really referring to her program for proving Fermat's Last Theorem.}
We know, too, that Germain wrote notes in 1822 on Libri's 1820
memoir,\footnote{Germain's three pages of notes \cite[cass. 7, ins.
56]{gernfl} \cite[p.~233]{delcentinabook}, while not directly about Fermat's
Last Theorem, do indicate an interest in modular solutions of roots of unity
equations, which is what encompasses the distribution of $p$-th powers modulo
$\theta$. Compare this with what she wrote to Gauss about Poinsot's work,
discussed in footnote \ref{poinsot}.} but this first memoir did not study
modular equations, hence was not relevant for the N-C condition. It seems
likely that she came to know of Libri's claims dooming her method, based
either on his presentations to the Academy in 1823/25 or the later memoir
published in 1829, particularly because Germain and Libri had met and were
personal friends from 1825 \cite[p. 117]{bucc} \cite[p. 140]{delcentinabook},
as well as frequent correspondents. It thus seems probable that sometime
between 1819 and 1825 Germain would have come to realize from Libri's work
that her grand plan could not work. However, we shall now see that she
determined this otherwise.

\subsubsection{Proof to Legendre that the plan fails for $p=3$}%

\begin{figure}
[pt]
\begin{center}
\includegraphics[
natheight=7.006700in,
natwidth=5.973300in,
height=5.9153in,
width=5.047in
]%
{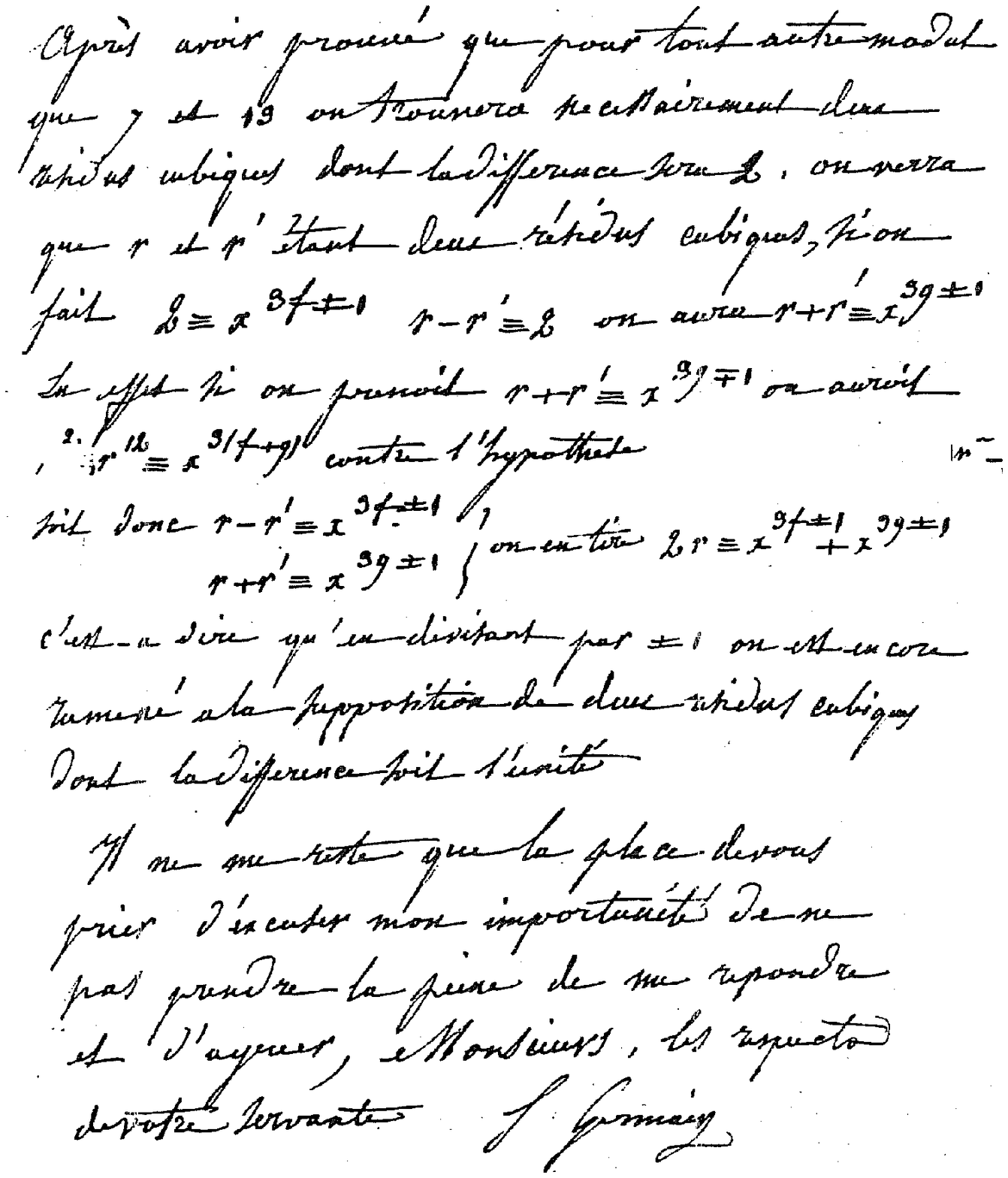}%
\caption{End of Germain's letter to Legendre}%
\label{legendre-l-image}%
\end{center}
\end{figure}

Beyond arguing as above that Germain very likely would have learned from
Libri's work that her grand plan cannot succeed, we have actually found
separate direct evidence of Germain's realization that her method of proving
Fermat's Last Theorem will not be successful, at least not in all cases.

While Manuscript A and her letter of 1819 to Gauss evince her belief that for
every prime $p>2$, there will be infinitely many auxiliary primes satisfying
the N-C condition, there is an undated letter to Legendre \cite{gerleg}
(described in the introduction) in which Germain actually proves the opposite
for $p=3$.

Sophie Germain began her three page letter by thanking Legendre for
\textquotedblleft telling\textquotedblright\ her \textquotedblleft
yesterday\textquotedblright\ that one can prove that all numbers of the form
$6a+1$ larger than $13$ have a pair of nonzero consecutive cubic residues.
This amounts to saying that for $p=3$, no auxiliary primes of the form
$\theta=2Np+1$ satisfy the N-C condition beyond $N=1,2$. At first sight this
claim is perplexing, since it seems to contradict Germain's success in
Manuscript A at proving Condition N-C for almost all odd primes $p$ whenever
$N=1,$ $2,$ $4,$ $5,$ $7,$ $8,$ $10$. However, the reader may check that for
$p=3$ her results in Manuscript A actually only apply for $N=1$ and $2$, once
one takes into account the exceptions, i.e., when either $\theta$ is not
prime, or Condition $2$-N-$p$ fails, or when she specifically excludes $p=3$
for $N=10$. So the claim by Legendre, mentioned in Germain's letter, that
there are only two valid auxiliary primes for $p=3$, is conceivably true.
Germain immediately writes a proof for him.

Since this proof is highly condensed, we will elucidate her argument here in
our own words, in modern terminology, and substantially expanded. Our aim is
to verify her claim, and at the same time experience the mathematical level
and sophistication of Germain's thinking. Figure \ref{legendre-l-image}
displays the end of the letter. The reader may notice that her last paragraph
of proof takes us fully twice as long to decipher and explain below.

\begin{GL}
For any prime $\theta$ of the form $6a+1,$ with $\theta>13$, there are
(nonzero) consecutive cubic residues. In other words, the N-C condition fails
for $\theta=2Np+1$ when $p=3$ and $N>2$, so the only valid auxiliary primes
for $p=3$ for the N-C condition are $\theta=7$ and $13.$
\end{GL}

\begin{pf}
We consider only the nonzero residues $1,\dots,6a$. Suppose that N-C is true,
i.e., there are no consecutive pairs of cubic residues (c.r.) amongst these,
and suppose further that there are also no pairs of c.r. whose difference is
$2$. (Note something important here. We mean literally residues, not
congruence classes, with this assumption, since obviously $1$ and $-1$ are
cubic congruence classes whose difference is $2$. But they are not both actual
residues, and their residues do not have difference $2$. So they do not
violate our assumption.) There are $2a$ c.r. distributed somehow amongst the
$6a$ residues, and without any differences of $1$ or $2$ allowed, according to
what we have assumed. Therefore to separate adequately these $2a$ residues
from each other there must be $2a-1$ gaps containing the $4a$ nonzero
non-cubic residues (n.c.r.), each gap containing at least $2$ n.c.r. Since
each of these $2a-1$ gaps has at least $2$ n.c.r., utilizing $4a-2$ n.c.r.,
this leaves flexibility for allocating only $2$ remaining of the $4a$ n.c.r.
This means that all the gaps must contain exactly $2$ n.c.r. except for either
a single gap with $4$ n.c.r., or two gaps with $3$ n.c.r. in each.

We already know of the specific c.r. $1$ and $8$ (recall $\theta=6a+1>13$).
and we know that $2$ and $3$ cannot be c.r. by our two assumptions. If $4$
were a c.r., then so would $8/4=2$ (alternatively, $8-4=4$ would violate N-C),
so $4$ is also not a c.r. Now Germain writes down a pattern for the sequence
of c.r. that we do not understand, and claims it is obviously absurd for
$\theta>13.$\footnote{Germain writes that the list is (presumably omitting
those at the ends) $1+4,$ $5+3,$ $8+3,$ $11+3,$ $14+3,\dots,6a-17,$ $6a-4$
[sic]$,$ $6a-11,$ $6a-8,$ $6a-5.$} We can easily arrive at a pattern and an
absurdity ourselves. From what Germain already has above, the c.r. sequence
must clearly be the list $1,$ $5,$ $8,$ $11,\dots,6a-10,$ $6a-7,$ $6a-4,$
$6a$, since the c.r. are symmetrically placed via negation modulo
$\theta=6a+1,$ and we know the gap sizes. Notice that the two exceptional gaps
must be of $3$ missing numbers each, located at the beginning and end. To see
this is absurd, consider first, for $\theta\geq6\cdot5+1=31$, the c.r.
$3^{3}=27$. Notice it contradicts the pattern listed above, since it is less
than $6a\geq30,$ but is not congruent to $2$ modulo $3,$ as are all the lesser
residues in the list except $1$. Finally, the only other prime $\theta>13$ is
$19,$ for which $4^{3}=64$ has residue $7$, which is not in the list.

So one of the two initial assumptions must be false. If N-C fails, we are
done. Therefore consider the failure of the other assumption, that there are
no pairs of c.r. whose difference is $2$. Let then $r$ and $r^{\prime}$ be
c.r. with $r-r^{\prime}=2$. Let $x$ be a primitive root of unity modulo
$\theta$, i.e., a generator of the cyclic group of multiplicative units
represented by the nonzero prime residues. We must have $2\equiv x^{3f\pm1},$
i.e., the power of $x$ representing $2$ cannot be divisible by $3$, since $2$
is not a c.r.

Now consider $r+r^{\prime}$. We claim that $r+r^{\prime}\not \equiv 0$, since
if $r+r^{\prime}\equiv0$, then $2=r-r^{\prime}\equiv r-\left(  -r\right)
=2r,$ yielding $r\equiv1$, and hence $r=1,$ which violates $r-r^{\prime}=2$.
Here it is critical to recall that we are dealing with actual residues $r$ and
$r^{\prime}$, both nonnegative numbers less than $6a+1$, i.e., the
requirements $r\equiv1$ and $r-r^{\prime}=2$ are incompatible, since there are
no $0<r,r^{\prime}<6a+1$ for which $r\equiv1$ and $r-r^{\prime}=2$; this is
related to the observation at the beginning that the congruence classes $1$
and $-1$ are not violating our initial assumption.

Since $r+r^{\prime}\not \equiv 0$, it is a unit, and thus must be congruent to
some power $x^{m}$. If $m$ were divisible by $3$, then the congruence
$r+r^{\prime}\equiv x^{m}$ would provide a difference of c.r. yielding another
c.r., which violates N-C after division by the latter. So we have
$r+r^{\prime}\equiv x^{3g\pm1}$. Now the sign in $3f\pm1$ must agree with that
in $3g\pm1$, since if not, say $r+r^{\prime}\equiv x^{3g\mp1}$, then
$r^{2}-r^{\prime2}=\left(  r-r^{\prime}\right)  \left(  r+r^{\prime}\right)
\equiv2x^{3g\mp1}\equiv x^{3f\pm1}x^{3g\mp1}=x^{3\left(  f+g\right)  }$, again
producing a difference of c.r. equal to another c.r., a contradiction.
Finally, we combine $r-r^{\prime}\equiv x^{3f\pm1}$ with $r+r^{\prime}\equiv
x^{3g\pm1}$ to obtain $2r\equiv x^{3f\pm1}+x^{3g\pm1}$, and thus $x^{3f\pm
1}r\equiv x^{3f\pm1}+x^{3g\pm1}$, becoming $r\equiv1+x^{3\left(  g-f\right)
}$, again contradicting N-C. Thus the original assumption of Condition N-C
must have been false. \quad\textsc{q.e.d.}
\end{pf}

This is quite impressive for a proof developed overnight.

These dramatic failures of Condition N-C for $p=3$ presumably greatly sobered
Germain's previous enthusiasm for pursuing her grand plan any further. We
mention in passing that, optimistic as Germain was at one point about finding
infinitely many auxiliary primes for each $p$, not only is that hope dashed in
her letter to Legendre, and by Libri's results, but even today it is not known
whether, for an arbitrary prime $p$, there is even one auxiliary prime
$\theta$ satisfying Condition N-C \cite[p.~301]{ribenboim-amateurs}.

\subsection{Germain's grand plan in other authors}

We know of no concrete evidence that anyone else ever pursued a plan similar
to Sophie Germain's for proving Fermat's Last Theorem, despite the fact that
Libri wrote of several (unnamed) mathematicians who attempted this method.
Germain's extensive work on this approach appears to be entirely,
independently, and solely hers, despite the fact that others were interested
in establishing Condition N-C for different purposes. In this section we will
see how and why other authors worked on Condition N-C, and compare with
Germain's methods.

\subsubsection{Legendre's methods for establishing Condition N-C}

Legendre did not mention Germain's full scale attack on Fermat's Last Theorem
via Condition N-C in his memoir of 1823, and we will discuss this later, when
we evaluate the interaction between Germain and Legendre in Section
\ref{Leg-Germ-interaction}. However, even ignoring any plan to prove Fermat's
Last Theorem outright, Legendre had two other reasons for wanting to establish
Condition N-C himself, and he develops N-C results in roughly the same range
for $N$ and $p$ as did Germain, albeit not mentioning her results.

One of his reasons was to verify Case 1 of Fermat's Last Theorem for many
prime exponents, since, recall, Condition N-C for a single auxiliary prime is
also one of the hypotheses of Sophie Germain's Theorem. Indeed, Legendre
develops results for N-C, and for the second hypothesis of her theorem, that
enable him to find a qualifying auxiliary prime for each odd exponent
$p\leq197$, which extends the scope of the table he implicitly attributed to
Germain. Legendre goes on to use his N-C results for a second purpose as well,
namely to show for a few small exponents that any solutions to the Fermat
equation would have to be very large indeed. We will discuss this additional
use of N-C in the next section.

Having said that Legendre obtained roughly similar N-C conclusions as Germain,
why do we claim that her approach to N-C verification is entirely independent?
This is because Germain's method of analyzing and proving the N-C condition,
explained in brief above, is utterly unlike Legendre's.\footnote{Del Centina
\cite[p. 370]{delcentina-germain-flt} also remarks on this.} We illustrate
this by quoting Legendre's explanation of why Condition N-C is always
satisfied for $N=2$, i.e., for $\theta=4p+1$. As we quote Legendre, we caution
that even his notation is very different; he uses $n$ for the prime exponent
that Germain, and we, call $p$. Legendre writes

\begin{quotation}
One can also prove that when one has $\theta=4n+1$, these two conditions are
also satisfied. In this case there are $4$ residues $r$ to deduce from the
equation $r^{4}-1=0$, which divides into two others $r^{2}-1=0$, $r^{2}+1=0$.
The second, from which one must deduce the number $\mu,$ is easy to
resolve\footnote{From earlier in the treatise, we know that $\mu$ here means a
primitive fourth root of unity, which will generate the four $n$-th powers.};
because one knows that in the case at hand $\theta$ may be put into the form
$a^{2}+b^{2}$, it suffices therefore to determine $\mu$ by the condition that
$a+b\mu$ is divisible by $\theta$; so that upon omitting multiples of $\theta
$, one can make $\mu^{2}=-1$, and the four values of $r$ become $r=\pm\left(
1,\mu\right)  $.

From this one sees that the condition $r^{\prime}=r+1$ can only be satisfied
in the case of $\mu=2$, so that one has $\theta=5$ and $n=1$, which is
excluded. ... \cite[\S 25]{legendre}
\end{quotation}

We largely leave it to the reader to understand Legendre's reasoning here. He
does not use the congruence idea or notation that Germain had adopted from
Gauss, he focuses his attention on the roots of unity from their defining
equation, he makes no use of the $2$-N-$p$ condition, but he is interested in
the consequences of the linear form $4n+1$ necessarily having a certain
quadratic form, although we do not see how it is germane to his argument. In
the next case, for $N=4$ and $\theta=8n+1$, he again focuses on the roots of
unity equation, and claims that this time the prime $8n+1$ must have the
quadratic form $a^{2}+2b^{2}$, which then enters intimately into an argument
related to a decomposition of the roots of unity equation. Clearly Legendre's
approach is completely unlike Germain's. Recall that Germain disposed of all
the cases $N=1,$ $2,$ $4,$ $5$ in one fell swoop with the first application of
her analysis of permuted placements of pairs of consecutive $p$-th powers,
whereas Legendre laboriously builds his analysis of $2N$-th roots of unity up
one value at a time from $N=1$. In short, Legendre focuses on the $p$-th
powers as $2N$-th roots of unity, one equation at a time, while Germain does
not, instead studying their permutations as $p$-th powers more generally for
what it indicates about their placement, and aiming for mathematical induction
on $N$.\footnote{Despite the apparently completely disjoint nature of the
treatments by Germain and Legendre of the N-C condition, it is quite curious
that their writings have a common mistake. The failure of N-C for $p=3$ when
$N=7$ is overlooked in Legendre's memoir, whereas in Germain's manuscript, as
we noted above, she explicitly calculated the failure of $2$-N-$p$ (and thus
of N-C) for this same combination, but then nonetheless mistakenly listed it
as valid for N-C in her table.}

\subsubsection{Dickson rediscovers permutation methods for Condition N-C}

Many later mathematicians worked to extend verification of the N-C condition
for larger values of $N$.\footnote{Legendre went to $N=8$ and Germain to
$N=10$, and actually to $N=11$ in another very much rougher manuscript draft
\cite[pp.\ 209r--214v, 216r--218v, 220r--226r]{ger9114}.} Their aim was to
prove Case 1 of Fermat's Last Theorem for more exponents by satisfying the
hypotheses of Sophie Germain's Theorem.

In particular, in 1908 L.\ E. Dickson published two papers
\cite{dickson1908,dickson1908-2} (also discussed in \cite[p. 763]{dickson})
extending the range of verification for Condition N-C to $N<74,$ and also $76$
and $128$ (each $N$ excepting certain values for $p$, of course), with which
he was able to apply Sophie Germain's theorem to prove Case 1 for all
$p<6$,$857$.

In light of the fact that Germain and Legendre had completely different
methods for verifying Condition N-C, one wonders what approach was taken by
Dickson. Dickson comments directly that his method for managing many cases
together has \textquotedblleft obvious advantages over the procedure of
Legendre\textquotedblright\ \cite[p. 27]{dickson1908-2}. It is then amazing to
see that his method is based directly (albeit presumably unbeknownst to him)
on the same theoretical observation made by Sophie Germain, that pairs of
consecutive $p$-th powers are permuted by two transformations of inversion and
subtraction to produce six more. He recognizes that these transformations form
a group of order six, which he calls the cross-ratio group (it consists of the
transformations of the cross-ratio of four numbers on the real projective line
obtained by permuting its variables \cite[pp.\ 112--113]{stillwell}), and is
isomorphic to the permutations on three symbols). Dickson observes that the
general form of these transformations of an arbitrary $p$-th power are the
roots of a sextic polynomial that must divide the roots of unity polynomial
for any $N$. This then forms the basis for much of his analysis, and even the
ad hoc portions have much the flavor of Germain's approach for $N>5$. In sum,
we see that Dickson's approach to the N-C condition more than three-quarters
of a century later could have been directly inspired by Germain's, had he
known of it.

\subsubsection{Modern approaches using Condition N-C}

Work on verifying the N-C condition has continued up to the close of the
twentieth century, largely with the aim of proving Case 1 using extensions of
Sophie Germain's Theorem.

By the middle of the 1980s results on the distribution of primes had been
combined with extensions of Germain's theorem to prove Case 1 of Fermat's Last
Theorem for infinitely many prime exponents \cite{adleman,fouvry}. It is also
remarkable that at least one yet more recent effort still harks back to what
we have seen in Germain's unpublished manuscripts. Recall that Germain
explained her intent to prove the N-C condition by induction on $N$. This is
precisely what a recent paper by David Ford and Vijay Jha does \cite{ford},
using some modern methods and computing power to prove by induction on $N$
that Case 1 of Fermat's Last Theorem holds for any odd prime exponent $p$ for
which there is a prime $\theta=2Np+1$ with $3\nmid N$ and $N\leq500$.

\subsection{Comparing Manuscripts A and D: Polishing for the prize
competition?}

We have analyzed Sophie Germain's grand plan to prove Fermat's Last Theorem,
which occupies most of Manuscript A. Manuscript D has the same title and
almost identical mathematical content and wording. Why did she write two
copies of the same thing? We can gain some insight into this by comparing the
two manuscripts more closely.

Manuscript D gives the impression of an almost finished exposition of
Germain's work on Fermat's Last Theorem, greatly polished in content and
wording over other much rougher versions amongst her papers. And it is
perfectly readable. However, it is not yet physically beautiful, since Germain
was clearly still refining her wording as she wrote it. In many places words
are crossed out and she continues with different wording, or words are
inserted between lines or in the margins to alter what has already been
written. There are also large parts of some pages left blank. By contrast,
Manuscript A appears essentially perfect. It is copied word for word almost
without exception from Manuscript D. It seems clear that Manuscript A was
written specifically to provide a visually perfected copy of Manuscript D.

One aspect of Manuscript D is quite curious. Recall that Manuscript A contains
a table with all the values for auxiliary primes satisfying Condition N-C for
$N\leq10$ and $3<p<100$. Germain explicitly introduces this table, referring
both ahead and back to it in the text, where it lies on page 17 of 20.
Manuscript D says all these same things about the table, but where the table
should be there is instead simply a side of a sheet left blank. Thus Germain
refers repeatedly to a table that is missing in what she wrote. This suggests
that as Germain was writing Manuscript D, she knew she would need to recopy it
to make it perfect, so she didn't bother writing out the table at the time,
saving the actual table for Manuscript A.

This comparison between Manuscripts A and D highlights the perfection of
presentation Sophie Germain sought in producing Manuscript A. Is it possible
that she was preparing this manuscript for submission to the French Academy
prize competition on the Fermat problem, which ran from 1816 to 1820? We will
discuss this further in Section \ref{prize-competition}.

\section{Large size of solutions\label{largesize}}

While Germain believed that her grand plan could prove Fermat's Last Theorem
for infinitely many prime exponents, she recognized that it had not yet done
so even for a single exponent. She thus wrote that she wished at least to show
for specific exponents that any possible solutions to the Fermat equation
would have to be extremely large.%

\begin{figure}
[pt]
\begin{center}
\includegraphics[
natheight=3.337300in,
natwidth=7.328400in,
height=2.3134in,
width=5.047in
]%
{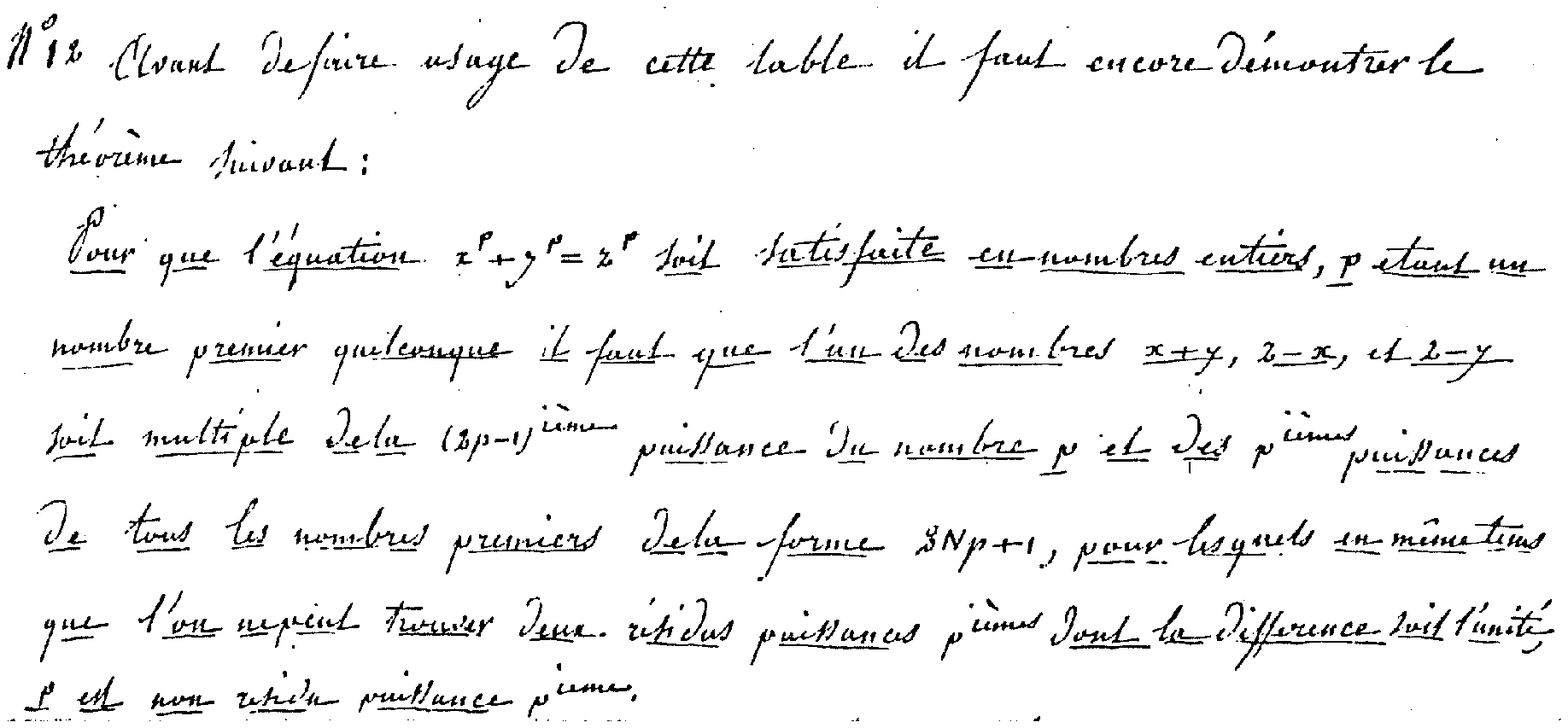}%
\caption{Beginning of the final section of Manuscript A, p. 17}%
\label{man-a-206-image}%
\end{center}
\end{figure}

In the last four pages of Manuscript A, Germain states, proves and applies a
theorem intended to accomplish this (Figure \ref{man-a-206-image}). She
actually states the theorem twice, first near the beginning of the manuscript
(Manuscript A, p. 3), where she recalls that any auxiliary prime satisfying
Condition N-C will have to divide one of the numbers $x$, $y$, $z$ in the
Fermat equation, but observes that to produce significant lower bounds on
solutions this way, one would need to employ rather large auxiliary primes.
Then she says

\begin{quotation}
\noindent fortunately one can avoid such impediment by means of the following
theorem:\footnote{\textquotedblleft heureusement on peut \'{e}viter un pareil
embarras au moyen du th\'{e}or\`{e}me suivant:\textquotedblright\ (Manuscript
A, p. 3)}
\end{quotation}

\begin{LS}
\textquotedblleft For the equation $x^{p}+y^{p}=z^{p}$ to be satisfied in
whole numbers, $p$ being any [odd] prime number, it is necessary that one of
the numbers $x+y$, $z-y$, and $z-x$ be a multiple of the ($2p-1$)$^{\text{th}%
}$ power of the number $p$ and of the $p^{\text{th}}$ powers of all the prime
numbers of the form [$\theta=$]$Np+1$, for which, at the same time, one cannot
find two $p^{\text{th}}$ power residues [mod $\theta$] whose difference is
one, and $p$ is not a $p^{\text{th}}$ power residue [mod $\theta
$].\textquotedblright\footnote{\textquotedblleft Pour que l'\'{e}quation
$x^{p}+y^{p}=z^{p}$ soit satisfaite en nombres entiers, $p$ \'{e}tant un
nombre premier quelconque; il faut que l'un des nombres $x+y$, $z-y$ et $z-x$
soit multiple de la ($2p-1$)$^{\text{i\`{e}me}}$ puissance du nombre $p$ et
des $p^{\text{i\`{e}mes}}$ puissances de tous les nombres premiers de la forme
$Np+1$, pour lesquels, en m\^{e}me tems que l'on ne peut trouver deux
r\'{e}sidus $p^{\text{i\`{e}mes}}$ puissances dont la difference soit
l'unit\'{e}, $p$ est non r\'{e}sidu puissance $p^{\text{i\`{e}me}}%
$.\textquotedblright\ (Manuscript A, p. 3 and p. 17)}
\end{LS}

\noindent(N.B: The theorem implicitly requires that at least one such $\theta$ exists.)

It is this theorem to which Germain was undoubtedly referring when, as we
noted earlier, she wrote to Gauss that any possible solutions would consist of
numbers \textquotedblleft whose size frightens the
imagination\textquotedblright. Early in Manuscript A she says that she will
apply the theorem for various values of $p$ using her table. She mentions here
that even just for $p=5$, the valid auxiliary primes $\theta=11,$ $41,$ $71,$
$101$ show that any solution to the Fermat equation would force a solution
number to have at least $39$ decimal digits.

We will see below that, as given, the proof of Germain's Large Size theorem is
insufficient, and we will discuss approaches she made to remedy this, as well
as an approach by Legendre to large size of solutions. But we will also see
that Sophie Germain's Theorem, the result she is actually known for today,
validly falls out of her proof.

\subsection{Germain's proof of large size of solutions\label{equations}}

Note first that the two hypotheses of Germain's Large Size theorem are the
same N-C condition she already studied at length for her grand plan, and a second:

\begin{pNp}
[$p$ is \textbf{N}ot a $p$-th power]$p$ is \textbf{n}ot a $p^{\text{th}}$
power residue, modulo $\theta$.
\end{pNp}

Of course this is precisely the second hypothesis of Sophie Germain's Theorem.

We now present a direct English translation of Germain's proof.

\subsubsection{The Barlow-Abel equations}

The proof implicitly begins with the fact that the N-C condition implies that
one of the numbers $x,y,z$ has to be divisible by $\theta$. We also provide
additional annotation, since Germain assumes the reader is already quite
familiar with many aspects of her equations.

\begin{quotation}
Assuming the existence of a single number subject to the double condition, I
will prove first that the particular number $x,y$ or $z$ in the equation
$x^{p}+y^{p}=z^{p}$ which is a multiple of the assumed number [$\theta$], must
necessarily also be a multiple of the number $p^{2}$.

Indeed, if the numbers $x,y,z$ are [assumed to be] relatively prime, then the
[pairs of] numbers
\[%
\begin{array}
[c]{lll}%
x+y & \mbox{and} & \ \ \ x^{p-1}-x^{p-2}y+x^{p-3}y^{2}-x^{p-4}y^{3}%
+\ \mbox{etc}\\
z-y & \mbox{and} & \ \ \ z^{p-1}+z^{p-2}y+z^{p-3}y^{2}+z^{p-4}y^{3}%
+\ \mbox{etc}\\
z-x & \mbox{and} & \ \ \ z^{p-1}+z^{p-2}x+z^{p-3}x^{2}+z^{p-4}x^{3}%
+\ \mbox{etc}.
\end{array}
\]
can have no common divisors other than $p$.\footnote{\textquotedblleft En
supposant l'existence d'un seul des nombres assujettis \`{a} cette double
condition, je prouverai d'abord que celui des nombres $x,$ $y$ et $z$ qui dans
l'\'{e}quation $x^{p}+y^{p}=z^{p}$ sera multiple du nombre suppos\'{e}, devra
necessairement \^{e}tre en m\^{e}me tems multiple du nombre $p^{2}$.
\par
\textquotedblleft En effet lorsque $x,$ $y$ et $z$ sont premiers entr'eux, les
nombres%
\[%
\begin{array}
[c]{lll}%
x+y & \mbox{et} & \ \ \ x^{p-1}-x^{p-2}y+x^{p-3}y^{2}-x^{p-4}y^{3}%
+\ \mbox{etc}\\
z-y & \mbox{et} & \ \ \ z^{p-1}+z^{p-2}y+z^{p-3}y^{2}+z^{p-4}y^{3}%
+\ \mbox{etc}\\
z-x & \mbox{et} & \ \ \ z^{p-1}+z^{p-2}x+z^{p-3}x^{2}+z^{p-4}x^{3}%
+\ \mbox{etc}.
\end{array}
\]
ne peuvent avoir d'autres diviseurs communs que le nombre $p$%
.\textquotedblright\ (Manuscript A, p. 18)}
\end{quotation}

For the first pair, this last claim can be seen as follows (and similarly for
the other pairs). Denote the right hand expression on the first line by
$\varphi(x,y)$. If some prime $q$ other than $p$ divides both numbers, then
$y\equiv-x\pmod q$, yielding $\varphi(x,y)\equiv px^{p-1}\pmod q$. Then $x$
and $x+y$ are both divisible by $q$, contradicting the assumption that $x$ and
$y$ are relatively prime. This excludes all primes other than $p$ as potential
common divisors of $x+y$ and $\varphi(x,y)$.

\begin{quotation}
If, therefore, the three numbers $x,y$, and $z$ were all prime to $p$, then
one would have, letting $z=lr,\ x=hn,\ y=vm$:\footnote{\noindent
\textquotedblleft Si on voulait donc que les trois nombres $x,y$, et $z$
fussent tous premiers a $p$ on aurait, en fesant $z=lr,$ $x=hn,$ $y=vm$:
\par%
\begin{align*}
x+y  &  =l^{p}\ \ \ \ \ \ \ x^{p-1}-x^{p-2}y+x^{p-3}y^{2}-x^{p-4}%
y^{3}+\ \mbox{etc}\ =r^{p}\\
z-y  &  =h^{p}\ \ \ \ \ \ \ z^{p-1}+z^{p-2}y+z^{p-3}y^{2}+z^{p-4}%
y^{3}+\ \mbox{etc}\ =n^{p}\\
z-x  &  =v^{p}\ \ \ \ \ \ \ \ z^{p-1}+z^{p-2}x+z^{p-3}x^{2}+z^{p-4}%
x^{3}+\ \mbox{etc}\ =m^{p}\text{.\textquotedblright}%
\end{align*}
\par
(Manuscript A, p. 18)}%
\begin{align}
x+y  &  =l^{p}\ \ \ \ \ \ \ x^{p-1}-x^{p-2}y+x^{p-3}y^{2}-x^{p-4}%
y^{3}+\ \mbox{etc}\ =r^{p}\tag{1}\\
z-y  &  =h^{p}\ \ \ \ \ \ \ z^{p-1}+z^{p-2}y+z^{p-3}y^{2}+z^{p-4}%
y^{3}+\ \mbox{etc}\ =n^{p}\tag{2}\\
z-x  &  =v^{p}\ \ \ \ \ \ \ \ z^{p-1}+z^{p-2}x+z^{p-3}x^{2}+z^{p-4}%
x^{3}+\ \mbox{etc}\ =m^{p}\text{.} \tag{3}%
\end{align}

\end{quotation}

Equations like these were given by Barlow around 1810, and stated apparently
independently by Abel in 1823 \cite[ch. III]{ribenboim-amateurs}.

One can derive these equations as follows. In the first line, the assumption
that $x,y,z$ are each relatively prime to $p$, along with the Fermat equation,
forces $x+y$ and $\varphi(x,y)$ to be relatively prime. Since the product of
$x+y$ and $\varphi(x,y)$ is equal to $z^{p}$, each of them must therefore be a
$p$th power, as she writes. The other lines have parallel proofs.

\subsubsection{Divisibility by $p$}

The next part of Germain's proof will provide a weak form of Sophie Germain's
Theorem, proving that one of $x,$ $y,$ $z$ must be divisible by $p$.

\begin{quotation}
Without loss of generality I assume that it is the number $z$ which is a
multiple of the prime number [$\theta$] of the form $2Np+1$, assumed to exist.
One therefore has that $l^{p}+h^{p}+v^{p}\equiv0\pmod
{2Np+1}$. And since by hypothesis there cannot be, for this modulus, two $p$th
power residues whose difference is $1$, it will be necessary that it is $l$
and not $r$, which has this modulus as a factor. Since $x+y\equiv
0\pmod{2Np+1}$, one concludes that $px^{p-1}\equiv r^{p}\pmod{2Np+1}$, that is
to say, because $x$ is a $p$th power residue, $p$ will also be a $p$th power
residue, contrary to hypothesis; thus the number $z$ must be a multiple of
$p$.\footnote{\textquotedblleft Pour fixer les id\'{e}es je supposerai que
c'est le nombre $z$ qui est multiple du nombre premier de la forme $2Np+1$
dont on a suppos\'{e} l'existence, on aura alors $l^{p}+h^{p}+v^{p}%
\equiv0\pmod
{2Np+1}$; et puisque par hypoth\`{e}se il ne peut y avoir pour ce module deux
r\'{e}sidus puissances $p^{\text{i\`{e}mes}}$ dont la difference soit
l'unit\'{e}, il faudra que ce soit $l$ et non par $r$ qui ait le m\^{e}me
module pour facteur. De $x+y\equiv0\pmod{2Np+1}$, on conclut $px^{p-1}\equiv
r^{p}\pmod{2Np+1}$ c'est \`{a} dire, \`{a} cause de $x$ r\'{e}sidu
$p^{\text{i\`{e}me}}$ puissance, $p$ aussi r\'{e}sidu $p^{\text{i\`{e}me}}$
puissance, ce qui est contraire \`{a} l'hypoth\`{e}se, il faut donc que le
nombre $z$ soit multiple de $p$.\textquotedblright\ (Manuscript A, p. 18)}
\end{quotation}

The N-C condition and the congruence $l^{p}+h^{p}+v^{p}\equiv0\pmod
{\theta=2Np+1}$ imply that either $l,\ h$, or $v$ is divisible by $\theta$. If
one of $h$ or $v$ were, then $x$ or $y$ would also be divisible by $\theta$,
contradicting the assumption that $x,\ y,\ z$ are relatively prime. This
implies that $l$ is the number divisible by $\theta$, and thus $y\equiv
-x\pmod{\theta}$. Substituting, we have $\varphi(x,y)\equiv px^{p-1}\equiv
r^{p}\pmod{\theta}$, as claimed. Furthermore, since $z\equiv0\pmod{\theta}$,
we conclude from $z-x=v^{p}$ that $x$ is a $p$th power modulo $\theta$.
Therefore, $p$ is also a $p$th power modulo $\theta$, a contradiction to the
other hypothesis of the theorem.

Thus we have derived a contradiction to the assumption that $x,\ y,\ z$ are
all prime to $p$, which indeed forces one of $x,\ y,\ z$ to be a multiple of
$p$. This is already the weak form of Sophie Germain's Theorem. But it is not
clear why $z$, the number divisible by $\theta$, has to be the one divisible
by $p$; this uncertainty is indicative of a flaw we will shortly observe.

In order to continue the proof, Germain now in effect implicitly changes the
assumption on $z$ to be that $z$ is the number known to be divisible by $p$,
but not necessarily by $\theta$, which in principle is fine, but must be kept
very clear by us. She replaces the first pair of equations by a new pair,
reflecting this change. (The remaining equations still hold, since $x$ and $y$
must be relatively prime to $p$.)

\subsubsection{Sophie Germain's Theorem as fallout}

Next in her proof comes the stronger form of Sophie Germain's Theorem.

\begin{quotation}
Setting actually $z=lrp$, the only admissible assumption is that
\begin{equation}
x+y=l^{p}p^{p-1},\qquad x^{p-1}-x^{p-2}y+x^{p-3}y^{2}-x^{p-4}y^{3}%
+\mbox{etc}=pr^{p}. \tag{1$^\prime$}%
\end{equation}
Because if, on the contrary, one were to assume that
\[
x+y=l^{p}p,\qquad x^{p-1}-x^{p-2}y+x^{p-3}y^{2}-x^{p-4}y^{3}%
+\mbox{etc}=p^{p-1}r^{p},
\]
then
\[
(x+y)^{p-1}-\{x^{p-1}-x^{p-2}y+x^{p-3}y^{2}+\mbox{etc}\}
\]
would be divisible by $p^{p-1}$. Observe that in the equation $2z-x-y=h^{p}%
+v^{p}$ the form of the right-hand side forces it to be divisible by $p$ or
$p^{2}$. Consequently, one sees that with the present assumptions $z$ has to
be a multiple of $p^{2}$.\footnote{\textquotedblleft En prenant actuellement
$z=lrp,$ la seule supposition admissible est
\[
x+y=l^{p}p^{p-1},\qquad x^{p-1}-x^{p-2}y+x^{p-3}y^{2}-x^{p-4}y^{3}%
+\mbox{etc}=pr^{p},
\]
car si on fesait au contraire%
\[
x+y=l^{p}p,\qquad x^{p-1}-x^{p-2}y+x^{p-3}y^{2}-x^{p-4}y^{3}%
+\mbox{etc}=p^{p-1}r^{p},
\]%
\[
(x+y)^{p-1}-\{x^{p-1}-x^{p-2}y+x^{p-3}y^{2}+\mbox{etc}\}
\]
serait divisible par $p^{p-1}$, parcons\'{e}quent si on observe que dans
l'\'{e}quation $2z-x-y=h^{p}+v^{p}$ la forme du second membre veut qu'il soit
premier a $p$, ou multiple de $p^{2}$ on verra que, dans les suppositions
presentes, $z$ aussi doit \^{e}tre multiple de $p^{2}$.\textquotedblright%
\ (Manuscript A, p. 18)}
\end{quotation}

To see Germain's first assertion one can argue as follows. Since $z^{p}%
=x^{p}+y^{p}$ must be divisible by $p$, we need only show that $\varphi(x,y)$
is divisible by exactly the first power of $p$. If we set $x+y=s$, then
\[
\varphi(x,y)=\frac{(s-x)^{p}+x^{p}}{s}=s^{p-1}-\binom{p}{1}s^{p-2}%
x+\cdots-\binom{p}{p-2}sx^{p-2}+\binom{p}{p-1}x^{p-1}.
\]
Now observe that all but the last summand of the right-hand side is divisible
by $p^{2}$, since $p$ divides $s=x+y\equiv x^{p}+y^{p}=z^{p}$ $\left(
\operatorname{mod}p\right)  $ by Fermat's Little Theorem, whereas the last
summand is divisible by exactly $p$, since $x$ is relatively prime to $p$.

Finally, to see that this forces $z$ to be divisible by $p^{2}$, observe that
the equation $2z-x-y=h^{p}+v^{p}$ ensures that $p$ divides $h^{p}+v^{p}$.
Furthermore, $p$ divides $h+v$ by Fermat's Little Theorem, applied to $h$ and
$v$. Now note that, since $h\equiv-v\pmod{p}$, it follows that $h^{p}%
\equiv-v^{p}\pmod{p^2}$. Thus $p^{2}$ divides $z$, since $p^{2}$ divides $x+y$
by Germain's new first pair of equations above.

This much of her proof constitutes a valid demonstration of what is called
Sophie Germain's Theorem.

\subsubsection{A mistake in the proof\label{mistake}}

Germain continues on to prove the further divisibility she claims by $\theta$.

\begin{quotation}
The only thing that remains to be proven is that all prime numbers of the form
$[\theta=]2Np+1$, which are subject to the same conditions as the number whose
existence has been assumed, are necessarily multiples [sic]\footnote{Germain
wrote \textquotedblleft multiples\textquotedblright\ here, but presumably
meant \textquotedblleft divisors\textquotedblright.} of $z$.

In order to obtain this let us suppose that it is $y$, for example, and not
$z$, that has one of the numbers in question as a factor. Then for this
modulus we will have $h^{p}-l^{p}\equiv v^{p}$, consequently $v\equiv
0,\ z\equiv x,\ pz^{p-1}\equiv m^{p}$, that is to say, $p$ is a $p$th power
residue contrary to the hypothesis.\footnote{\textquotedblleft La seule chose
qui reste \`{a} prouver est que tous les nombres premier de la forme $2Np+1$
qui sont assujettis aux m\^{e}mes conditions que celui de la m\^{e}me forme
dont en a suppos\'{e} l'existence sont necessairement multiples [sic] de $z$.
\par
\textquotedblleft Pour y parvenir supposons que ce soit $y$, par exemple et
non pas $z$, qui ait un des nombres dont il s'agit pour facteur, nous aurons
pour ce module $h^{p}-l^{p}\equiv v^{p}$, parcons\'{e}quent $v\equiv
0,\ z\equiv x,\ pz^{p-1}\equiv m^{p}$, c'est a dire $p$ residu puissance
$p^{\text{i\`{e}me}}$ contre l'hypoth\`{e}se.\textquotedblright\ (Manuscript
A, pp. 18--19)}
\end{quotation}

Here Germain makes a puzzling mistake.\footnote{Del Centina \cite[p.
365ff]{delcentina-germain-flt} does not seem to notice this mistake.} Rather
than using the equation (1$^{\prime}$), resulting from the $p$-divisibility
assumption on $z$, she erroneously uses the original equation (1) which
required the assumption that all of $x,y,z$ are relatively prime to $p$.
Subtracting (1) from (2) and comparing the result to (3), she obtains the
congruence $h^{p}-l^{p}\equiv v^{p}\pmod{\theta}$, since $y\equiv
0\pmod{\theta}$. Although this congruence has been incorrectly obtained, we
will follow how she deduces from it the desired contradiction, partly because
we wish to see how the entire argument might be corrected. Since neither $h$
nor $l$ can be divisible by $\theta$ (since neither $x$ nor $z$ are), the N-C
Condition implies that $v\equiv0\pmod{\theta}$, hence $z\equiv x$. Thus,
$pz^{p-1}\equiv m^{p}$ follows from the right-hand equation of (3). Further,
$z\equiv h^{p}$ follows from (2), since $y\equiv0$, and, finally, this allows
the expression of $p$ as the residue of a $p$-th power, which contradicts the
$p$-N-$p$ Condition.

Except for the mistake noted, the proof of Germain's theorem is complete. If
instead the correct new equation (1$^{\prime}$) had been used, then in place
of the N-C Condition, the argument as written would need a condition analogous
to N-C, but different, for the congruence
\[
h^{p}-l^{p}p^{p-1}\equiv v^{p}%
\]
resulting from subtracting (1$^{\prime}$) from (2) instead of (1) from (2).
That is, we could require the following additional hypothesis:

\begin{Np-inv}
[\textbf{N}o $p^{-1}$ differences]There are \textbf{n}o two nonzero
$p^{\mathrm{th}}$-power residues that differ by $p^{-1}$ (equivalently, by
$-2N$) modulo $\theta$.
\end{Np-inv}

Clearly, adding this condition as an additional hypothesis would make the
proof of the theorem valid.

\subsubsection{Attempted remedy}

Did Germain ever realize this problem, and attempt to correct it?

To the left of the very well defined manuscript margin, at the beginning of
the paragraph containing the error, are written two words in much smaller
letters and a thicker pen. These words are either \textquotedblleft voyez
errata\textquotedblright\ or \textquotedblleft voyez erratu\textquotedblright.
This is one of only four places in Manuscript A where marginal notes mar its
visual perfection. None of these appears in Manuscript D, from which
Manuscript A was meticulously copied. So Germain saw the error in Manuscript
A, but probably later, and wrote an erratum about it. Where is the erratum?

Most remarkably, not far away in the same archive of her papers, tucked
apparently randomly in between other pages, we find two sheets \cite[pp.~214r,
215v]{ger9114} clearly titled \textquotedblleft errata\textquotedblright\ or
\textquotedblleft erratu\textquotedblright\ in the same writing style as the
marginal comment.

The moment one starts reading these sheets, it is clear that they address
precisely the error Germain made. After writing the corrected equations
(1$^{\prime}$), (2), (3) (in fact she refines them even more, incorporating
the $p^{2}$ divisibility she just correctly deduced) Germain notes that it is
therefore a congruence of the altered form
\[
l^{p}p^{2p-1}+h^{p}+v^{p}\equiv0
\]
that should hopefully lead to a contradiction. It is not hard to see that the
N-$p^{-1}$ and $p$-N-$p$ conditions will suffice for this, but Germain
observes right away that a congruence nullifying the N-$p^{-1}$ condition in
fact exists for the very simplest case of interest to her, namely $p=5$ and
$N=1$, since $1$ and $-1$ are both $5$-th powers, and they differ by
$2N=2$.\footnote{In fact the reader may check in various examples for small
numbers that the N-$p^{-1}$ condition seems to hold rather infrequently
compared with the N-C condition, so simply assuming the N-$p^{-1}$ condition
as a hypothesis makes a true theorem, but perhaps not a very useful one.}

Germain then embarks on an effort to prove her claim by other means, not
relying on assuming the N-$p^{-1}$ condition. She develops arguments and
claims based on knowledge of quadratic forms and quadratic reciprocity,
including marginal comments that are difficult to interpret. There is more
work to be done understanding her mathematical approach in this erratum, which
ends inconclusively. What Germain displays, though, is her versatility, in
bringing in quadratic forms and quadratic reciprocity to try to resolve the issue.

\subsubsection{Verifying Condition $p$-N-$p$: a theoretical approach}

We return now from Germain's erratum to discuss the end of Manuscript A.
Germain follows her Large Size of Solutions theorem with a method for finding
auxiliary primes $\theta$ of the form $2Np+1$ satisfying the two conditions
(N-C and $p$-N-$p$) required for applying the theorem.

Even though we now realize that her applications of the Large Size theorem are
unjustified, since she did not succeed in providing a correct proof of the
theorem, we will describe her methods for verifying its hypotheses, in order
to show their skill, their application to Sophie Germain's theorem, and to
compare them with the work of others.

Earlier in the manuscript Germain has already shown her methods for verifying
Condition N-C for her grand plan. She now focuses on verifying Condition
$p$-N-$p$, with application in the same range as before, i.e., for auxiliary
primes $\theta=2Np+1$ using relevant values of $N\leq10$ and odd primes
$p<100$.

Germain first points out that since $\theta=2Np+1$, therefore $p$ will be a
$p$-th power modulo $\theta$ if and only if $2N$ is also, and thus, due to the
cyclic nature of the multiplicative units modulo $\theta$, precisely if
$\left(  2N\right)  ^{2N}-1$ is divisible by $\theta$. Yet before doing any
calculations of this sort, she obviates much effort by stating another
theoretical result: For $N$ of the form $2^{a}p^{b}$ in which $a+1$ and $b+1$
are prime to $p$, she claims that $p$ cannot be a $p$-th power modulo $\theta$
provided $2$ is not a $p$-th power modulo $\theta$. Of course the latter is a
condition ($2$-N-$p$) she already studied in detail earlier for use in her N-C
analyses. Indeed the claim follows because $2^{a+1}p^{b+1}=2Np\equiv\left(
-1\right)  ^{p}$, which shows that $2$ and $p$ must be $p$-th powers together
(although the hypothesis on $b$ is not necessary for just the implication she
wishes to conclude). Germain points out that this result immediately covers
$N=1,2,4,8$ for all $p.$ In fact, there is in these cases no need for Germain
even to check the $2$-N-$p$ condition, since she already earlier verified N-C
for these values of $N$, and $2$-N-$p$ follows from N-C. Germain easily
continues to analyze $N=5,7,10$ for Condition $p$-N-$p$ by factoring $\left(
2N\right)  ^{2N}-1$ and looking for prime factors of the form $2Np+1$.
Astonishingly, by this method Germain deduces that there is not a single
failure of Condition $p$-N-$p$ for the auxiliary primes $\theta=2Np+1$ in her
entire previously drawn table of values satisfying Condition N-C.

Germain ends Manuscript A by drawing conclusions on the minimum size of
solutions to Fermat equations for $2<p<100$ using the values for $\theta$ in
her table. Almost the most modest is her conclusion for $p=5$. Since her
techniques have verified that the auxiliaries $11,$ $41,$ $71,$ $101$ all
satisfy both Conditions N-C and $p$-N-$p$, Germain's Large Size theorem (if it
were true) ensures that if $x^{5}+y^{5}=z^{5}$ were true in positive numbers,
then one of the numbers $x+y$, $z-y$, $z-x$ must be divisible by $5^{9}%
11^{5}41^{5}71^{5}101^{5}$, which Germain notes has at least 39 decimal digits.

\subsection{Condition $p$-N-$p$ and large size in other authors}

Legendre's footnote credits Germain for Sophie Germain's Theorem and for
applying it to prove Case 1 for odd primes $p<100$ \cite[\S 22]{legendre}. For
the application he exhibits a table providing, for each $p$, a single
auxiliary prime satisfying both conditions N-C and $p$-N-$p$, based on
examination of a raw numerical listing of all its $p$-th power residues.

Thus he leaves the impression that Germain verified that her theorem was
applicable for each $p<100$ by brute force residue computation with a single
auxiliary. In fact, there is even such a residue table to be found in
Germain's papers \cite[p.~151v]{ger9114}, that gives lists of $p$-th power
residues closely matching Legendre's table.\footnote{There are a couple of
small differences between Legendre's table of residues and the one we find in
Germain's papers. Germain states that she will not list the residues in the
cases when $N\leq2$ in the auxiliary prime, suggesting that she already knew
that such auxiliary primes are always valid. And while Germain, like Legendre,
generally lists for each $p$ the residues for only the single smallest
auxiliary prime valid for both N-C and $p$-N-$p$, in the case of $p=5$ she
lists the residues for several of the auxiliaries that she validated in
Manuscript A.} Legendre's table could thus easily have been made from hers.
This, however, is not the full story, contrary to the impression received from Legendre.

\subsubsection{Approaches to Condition $p$-N-$p$}

Both Legendre and Germain analyze theoretically the validity of Condition
$p$-N-$p$ as well as that of N-C for a range of values of $N$ and $p$, even
though, as with Germain's grand plan for proving Fermat's Last Theorem via
Condition N-C, Legendre never indicates her efforts at proving large size for
solutions by finding multiple auxiliary primes satisfying both Conditions N-C
and $p$-N-$p$.

Moreover, since all Legendre's work at verifying N-C and $p$-N-$p$ comes after
his footnote crediting Germain, he is mute about Germain developing techniques
for verifying either condition. Rather, the clear impression his treatise
leaves to the reader is that Sophie Germain's Theorem and the brute force
table are hers, while all the techniques for verifying Conditions N-C and
$p$-N-$p$ are his alone.

As we have seen, though, Germain qualifies auxiliaries to satisfy both N-C and
$p$-N-$p$ entirely by theoretical analyses, and her table in Manuscript A has
no brute force listing of residues. In fact she developed general techniques
for everything, with very little brute force computation evident, and was very
interested in verifying her conditions for many combinations of $N$ and $p$,
not just one auxiliary for each $p$. In short, the nature of Legendre's credit
to Germain for proving Case 1 for $p<100$ leaves totally invisible her much
broader theoretical work that we have uncovered in Manuscript A.

We should therefore investigate, as we did earlier for Condition N-C, how
Legendre's attempts at verifying Condition $p$-N-$p$ compare with Germain's,
to see if they are independent.

\subsubsection{Legendre on Condition $p$-N-$p$}

Legendre's approach to verifying Condition $p$-N-$p$ for successive values of
$N$ is at first rather ad hoc, then based on the criterion whether $\theta$
divides $p^{2N}-1,$ slowly evolving to the equivalent divisibility of $\left(
2N\right)  ^{2N}-1$ instead, and appeals to his \emph{Th\'{e}orie des Nombres}
for\emph{ }finding divisors of numbers of certain forms.

Unlike Germain's methods, there is no recognition that many $N$ of the form
$2^{a}p^{b}$ are amenable to appeal to Condition $2$-N-$p$. Suffice it to say
that, as for Condition N-C, Legendre's approaches and Germain's take different
tacks, with Germain starting with theoretical transformations that make
verification easier, even though in the end they both verify Condition
$p$-N-$p$ for roughly the same ranges of $N$ and $p$. There are aspects with
both the N-C and $p$-N-$p$ analyses where Germain goes further than Legendre
with values of $N$ and $p$, and vice versa.

Even their choices of symbols and notation are utterly different. Legendre
never uses the congruence notation that Gauss had introduced almost a quarter
century before, while Germain is fluent with it. Legendre quotes and relies on
various results and viewpoints from the second edition of his
\emph{Th\'{e}orie des Nombres}, and never considers Condition $2$-N-$p$ either
for N-C or $p$-N-$p$ analysis, whereas it forms a linchpin in Germain's
approach to both. Germain rarely refers to Legendre's book or its results, but
uses instead her intimate understanding of the multiplicative structure of
prime residues from Gauss's \emph{Disquisitiones}.

We are left surprised and perplexed by the lack of overlap in mathematical
approach between Germain's Manuscript A and Legendre's treatise, even though
the two are coming to the same conclusions page after page. There is nothing
in the two manuscripts that would make one think they had communicated, except
Legendre's footnote crediting Germain with the theorem that today bears her
name. It is as though Legendre never saw Germain's Manuscript A, a thought we
shall return to below. Four factors leave us greatly perplexed at this
disparity. First, years earlier Legendre had given Germain his strong
mentorship during the work on elasticity theory that earned her a prize of the
French Academy. Second, Legendre's own research on Fermat's Last Theorem was
contemporaneous with Germain's. Third, Germain's letter to Legendre about the
failure of N-C for $p=3$ demonstrates detailed interaction. Fourth, we shall
discuss later that Legendre's credit to Germain does match quite well with her
Manuscript B. How could they not have been in close contact and sharing their
results and methods? In the end, at the very least we can conclude that each
did much independent work, and should receive separate credit for all the
differing techniques they developed for analyzing and verifying the N-C and
$p$-N-$p$ conditions.

\subsubsection{Legendre's approach to large size of solutions}

Legendre describes not just Sophie Germain's Theorem and applications, but
also large size results similar to Germain's, although he makes no mention of
his large size results having anything to do with her. Thus we should compare
their large size work as well.

Germain presents a theorem about large size, and quite dramatic specific
consequences, but the theorem is flawed and her attempts at general repair
appear inconclusive. Legendre, like Germain, studies whether all qualifying
auxiliary primes $\theta$ must divide the same one of $x,$ $y,$ $z$ that
$p^{2}$ does, which is where Germain went wrong in her original manuscript.
Like Germain in her erratum, Legendre recognizes that the N-$p^{-1}$ condition
would ensure the desired $\theta$ divisibility. He, like Germain, also presses
on in an alternative direction, since the condition is not necessarily (in
fact perhaps not even often) satisfied. But here, just as much as in his
differing approach to verifying the N-C and $p$-N-$p$ conditions, Legendre
again chooses a completely different alternative approach than does Germain.

Legendre analyzes the placement of the $p$-th power residues more deeply in
relation to the various expressions in equations (1$^{\prime}$), (2), (3)
above, and finds additional conditions, more delicate than that of N-$p^{-1}$,
which will ensure the desired $\theta$ divisibility for concluding large size
of solutions. Specifically, for example, when $p=5$ Legendre has the same
auxiliaries $\theta=11,$ $41,$ $71,$ $101$ satisfying N-C and $p$-N-$p$ as had
Germain.\footnote{\label{legendre-to-1000}Although Legendre never mentions the
grand plan for proving Fermat's Last Theorem, he is interested in how many
valid auxiliaries there may be for a given exponent. He claims that between
$101$ and $1000$ there are no auxiliaries for $p=5$ satisfying the two
conditions, and that this must lead one to expect that $101$ is the last. This
presages Libri's claims that for each $p$ there are only finitely many
auxiliaries satisfying N-C, and is the one hint we find in Legendre of a
possible interest in the grand plan.} However, as Germain explicitly pointed
out for $\theta=11$ in her erratum, Condition N-$p^{-1}$ fails; in fact
Legendre's calculations show that it fails for all four auxiliaries. While
Germain attempted a general fix of her large size theorem using quadratic
forms and quadratic reciprocity, Legendre's delicate analysis of the placement
of $5$-th powers shows that $11,$ $71,$ $101$ (but not $41$) must divide the
same one of $x,$ $y,$ $z$ as $p^{2}$, and so he deduces that some sum or
difference of two of the indeterminates must be divisible by $5^{9}%
11^{5}71^{5}101^{5}$, i.e., must have at least $31$ digits. This is weaker
than the even larger size Germain incorrectly deduced, but it is at least a
validly supported conclusion. Legendre successfully carries this type of
analysis on to exponents $p=7,$ $11$, $13$, concluding that this provides
strong numerical evidence for Fermat's Last Theorem. But he does not attempt a
general theorem about large size of solutions, as did Germain. As with their
work on Conditions N-C and $p$-N-$p$, we are struck by the disjoint approaches
to large size of solutions taken by Germain and Legendre. It seems clear that
they each worked largely independently, and there is no evidence in their
manuscripts that they influenced each other.

\subsubsection{Rediscovery of Germain's approach to Condition $p$-N-$p$}

Later mathematicians were as unaware of Germain's theoretical analysis of
Condition $p$-N-$p$ as they were of her approach to Condition N-C, again
because Legendre's published approach was very different and introduced
nothing systematically helpful beyond basic calculation, and Germain's work
was never published \cite[ch. 8]{bucc}.

In particular, the fact that for values of $N$ of the form $2^{a}p^{b}$ for
which $p$ and $a$ are relatively prime, Condition $p$-N-$p$ follows from
$2$-N-$p$, was essentially (re)discovered by Wendt in 1894 \cite{wendt}, and
elaborated by Dickson \cite{dickson1908} and Vandiver\footnote{For
comprehensive views of Vandiver's contributions, especially in relation to
Case 1, see \cite{corry1,corry2}.} \cite{vandiver} in the twentieth century.

\section{Exponents of form $2(8n\pm3)$\label{specialform}}

We will consider now what we call Manuscript B, entitled
\emph{D\'{e}monstration de l'impossibilit\'{e} de satisfaire en nombres
entiers \`{a} l'\'{e}quation }$z^{2(8n\pm3)}=y^{2(8n\pm3)}+x^{2(8n\pm3)}$. By
the end of the manuscript, although it is written in a less polished fashion,
it is clear that Germain has apparently proven Fermat's Last Theorem for all
exponents of the form $2\left(  8n\pm3\right)  $, where $p=$ $8n\pm3$ is prime.

Germain states and proves three theorems, and then has a final argument
leading to the title claim. We shall analyze this manuscript for its approach,
for its connection to her other manuscripts and to Legendre's attribution to
her, and for its correctness.

Although Germain does not spell out the big picture, leaving the reader to put
it all together, it is clear that she is proceeding to prove Fermat's Last
Theorem via the division we make today, between Case 1 and Case 2, separately
eliminating solutions in which the prime exponent $p=$ $8n\pm3$ either does
not or does divide one of $x^{2},$ $y^{2}$, $z^{2}$ in the Fermat equation
$\left(  x^{2}\right)  ^{p}+\left(  y^{2}\right)  ^{p}=\left(  z^{2}\right)
^{p}$.

\subsection{Case 1 and Sophie Germain's Theorem}

Germain begins by claiming to eliminate solutions in which none are divisible
by $p$, and actually claims this for all odd prime exponents, writing

\begin{quotation}
First Theorem. For any [odd] prime number $p$ in the equation $z^{p}%
=x^{p}+y^{p}$, one of the three numbers $z$, $x$, or $y$ will be a multiple of
$p^{2}$.\footnote{\noindent\textquotedblleft Th\'{e}or\`{e}me premier.
\emph{Quelque soit le nombre premier }$p$\emph{ dans l'\'{e}quation }%
$z^{p}=x^{p}+y^{p}$\emph{ l'un des trois nombres }$z,$ $x$ $\emph{ou}$ $y$
\emph{sera multiple de }$p^{2}$.\textquotedblright\ (Manuscript B, p. 92r)}
\end{quotation}

Today we name this Case 1 of Fermat's Last Theorem, that solutions must be
$p$-divisible (Germain claims a little more, namely $p^{2}$ divisibility).
Note that there are no hypotheses as stated, since Germain wishes to evince
that Case 1 is true in general, and move on to Case 2 for the exponents at
hand. She does, however, immediately recognize that to prove this, she
requires something else:

\begin{quotation}
To demonstrate this theorem it suffices to suppose that there exists at least
one prime number $\theta$ of the form $2Np+1$ for which at the same time one
cannot find two $p^{\text{th}}$ power residues [mod $\theta$] whose difference
is one, and $p$ is not a $p^{\text{th}}$ power residue [mod $\theta
$].\footnote{\textquotedblleft Pour d\'{e}montrer ce th\'{e}or\`{e}me il
suffit de supposer qu'il existe au moins un nombre premier $\theta$ de la form
$2Np+1$ pour lequel en m\^{e}me tems que l'on ne peut trouver deux residus
puissances $p^{\text{i\`{e}me}}$ dont la difference soit l'unit\'{e} $p$ est
non residu puissance $p^{\text{i\`{e}me}}$.\textquotedblright\ (Manuscript B,
p. 92r)}
\end{quotation}

Today we recognize this as the hypothesis of Sophie Germain's Theorem, whereas
for her it was not just a hypothesis, but something she believed was true and
provable by her methods, since she goes on to say

\begin{quotation}
Not only does there always exist a number $\theta$ satisfying these two
conditions, but the course of calculation indicates that there must be an
infinite number of them. For example, if $p=5$, then $\theta=2\cdot5+1=11$,
\quad$2\cdot4\cdot5+1=41$, \quad$2\cdot7\cdot5+1=71$, \quad$2\cdot
10\cdot5+1=101$, etc.\footnote{\textquotedblleft Non seulement il existe
toujours un nombre $\theta$ qui satisfait \`{a} cette double condition mais la
marche du calcul indique qu'il doit s'entrouver une infinit\'{e} \quad
$p=5$\quad$\theta=2\cdot5+1=11$\textsf{, }\quad$2\cdot4\cdot5+1=41$\textsf{,
}\quad$2\cdot7\cdot5+1=71$\textsf{, }\quad$2\cdot10\cdot5+1=101$,
etc.\textquotedblright\ \ (Manuscript B, p. 92r}

\end{quotation}

Recall that Germain spends most of Manuscript A developing powerful techniques
that support this belief in Conditions N-C and $p$-N-$p$, and that confirm
them for $p<100$, so it is not surprising that she wishes to claim to have
proven Case 1 of Fermat's Last Theorem, even though she still recognizes that
there are implicit hypotheses she has not completely verified for all exponents.

Germain's proof of her First Theorem is much like the beginning of her proof
of the Large Size theorem of Manuscript A, which we laid out in Section
\ref{largesize}. Recall that the Large Size proof went awry only after the
$p^{2}$ divisibility had been proven, so her proof here,\footnote{The proof of
Theorem 1 in Manuscript B is largely reproduced, in translation, in \cite[p.
189ff]{lp}.} as there, proves $p^{2}$ divisibility without question. This is
the closest to an independent statement and proof we find in her manuscripts
of what today is called Sophie Germain's Theorem.

However, most curiously, at the end of the proof of the First Theorem she
claims also that the $p^{2}$ divisibility applies to the same one of $x,$ $y,$
$z$ that is divisible by the auxiliary prime $\theta$, which is the same as
the claim, ultimately inadequately supported, where her Large Size proof in
Manuscript A began to go wrong. While she makes no use of this additional
claim here (so that it is harmless to her line of future argument in this
manuscript), it leads us to doubt a conjecture one could otherwise make about
Manuscript B. One could imagine that the First Theorem was written down as a
means of salvaging what she could from the Large Size theorem, once she
discovered the flaw in the latter part of its proof. But since the confusion
linked to the flawed claim there appears also here (without proof), even
though without consequent maleffect, we cannot argue that this manuscript
contains a corrected more limited version of the Large Size theorem argument.

\subsection{Case 2 for $p$ dividing $z$}

The rest of Manuscript B deals with Case 2 of Fermat's Last Theorem, which is
characterized by equations (1$^{\prime}$), (2), (3) in Section \ref{equations}%
. For completeness, we mention that Theorem 2 contains a technical result not
relevant to the line of proof Germain is developing. Perhaps she placed it and
its proof here simply because it was a result of hers about Case 2, which is
the focus of the rest of the manuscript.\footnote{Theorem 2 asserts that in
the equations (1$^{\prime}$), (2), (3) pertaining in Case 2, the numbers $r,$
$m,$ $n$ can have prime divisors only of the form $2Np+1$, and that moreover,
the prime divisors of $r$ must be of the even more restricted form $2Np^{2}%
+1$. Legendre also credits this result to Germain in his footnote.}

As we continue with Case 2, notice that, by involving squares, the equation
$\left(  x^{2}\right)  ^{p}+\left(  y^{2}\right)  ^{p}=\left(  z^{2}\right)
^{p}$ has an asymmetry forcing separate consideration of $z$ from $x$ or $y$
in proving Fermat's Last Theorem. Germain addresses the first of these, the
$p$-divisibility of $z$, in her Theorem 3, which asserts that $z$ cannot be a
multiple of $p$, if $p$ has the form $8n+3,$ $8n+5,$ or $8n+7$. She proves
Theorem 3 by contradiction, by assuming that $z$ is divisible by $p$. Her
proof actually begins with some equations that require some advance
derivation. Using the relative primality of the key numbers in each pair of
the Case 2 equations (1$^{\prime}$), (2), (3) of Manuscript A, for pairwise
relatively prime solutions $x^{2},$ $y^{2},$ $z^{2}$ (once the extra $p^{2}$
divisibility is built in), the reader may easily verify that the left trio of
these equations becomes\footnote{We do not see how she obtains $4p-1$ as
exponent, rather than just $2p-1$, even after including the stronger $p^{2}$
divisibility; but $2p-1$ suffices.}%

\begin{align*}
x^{2}+y^{2}  &  =p^{4p-1}l^{2p}\\
z^{2}-y^{2}  &  =h^{2p}\\
z^{2}-x^{2}  &  =v^{2p}.
\end{align*}
The text of Germain's proof begins with these equations.

Germain quickly confirms Theorem 3 for $p=8n+3$ and $8n+7$ using the fact,
long known from Fermat's time, that a sum of squares can contain no prime
divisors of these two forms. For $p=8n+5$ she must argue differently, as follows.

Because $z-y$ and $z+y$ (respectively $z-x$ and $z+x$) are relatively prime,
one has $z+y=\left(  h^{\prime}\right)  ^{2p}$ and $z+x=\left(  v^{\prime
}\right)  ^{2p}$, whence $y^{2}\equiv\left(  h^{\prime}\right)  ^{4p}$ (mod
$p$) and $x^{2}\equiv\left(  v^{\prime}\right)  ^{4p}$ (mod $p$), yielding
$\left(  h^{\prime}\right)  ^{4p}+\left(  v^{\prime}\right)  ^{4p}\equiv0$
(mod $p$) since $x^{2}+y^{2}$ is divisible by $p$. This, she points out, is a
contradiction, since $-1$ is not a biquadratic residue modulo $8n+5$.

The unfortunate flaw in this proof is perhaps not obvious at first. The
$2p$-th power expressions for $z+y$ and $z+x$ rely on $z-y$ and $z+y$
(respectively $z-x$ and $z+x$) being relatively prime. This would be true from
the pairwise relative primality of $x,$ $y,$ $z,$ if the numbers in each
difference had opposite parity, but otherwise their difference and sum have
precisely $2$ as greatest common divisor. Writing $\left(  x^{p}\right)
^{2}+\left(  y^{p}\right)  ^{2}=\left(  z^{p}\right)  ^{2}$ and recalling
basics of Pythagorean triples, we see that opposite parity fails either for
$z-y$ or $z-x$. Suppose without loss of generality that it is $z-y$. Then
either $z-y$ or $z+y$ has only a single $2$ as factor (since $y$ and $z$ are
relatively prime), so it cannot be a $2p$-th power. One can include this
single factor of $2$ and redo Germain's analysis to the end, but one then
finds that it comes down to whether or not $-4$ is a biquadratic residue
modulo $8n+5$, and this unfortunately is true, rather than false as for $-1$.
So Germain's proof of Theorem 3 appears fatally flawed for $p=8n+5$.

\subsection{Case 2 for $p$ dividing $x$ or $y$}

In her final argument after Theorem 3, Germain finishes Case 2 for $p=8n+3$
and $8n-3$ by dealing with the second possible situation, where either $x$ or
$y$ is divisible by $p$. This argument again builds from enhanced versions of
equations similar to (1$^{\prime}$), (2), (3), but is considerably more
elaborate, rising up through detailed study of the specific cases $p=5,$ $13,$
$29,$ until she is able to end with an argument applying to all $p=8n+3$ and
$8n-3$. However, since the argument proceeds initially as did the proof of
Theorem 3, it too relies on the same mistaken assumption about relative
primality that misses an extra factor of $2$, and one finds that accounting
for this removes the contradiction Germain aims for, no matter what value $p$ has.

\subsection{Manuscript B as source for Legendre?}

In the end we must conclude that this proof of the bold claim to have proven
Fermat's Last Theorem for many exponents fails due to an elementary mistake.
But what is correct in Manuscript B fits extremely well with what Legendre
wrote about Germain's work. The manuscript contains precisely the correct
results Legendre credits to Germain, namely Sophie Germain's Theorem and the
technical result of Theorem 2 about the equations in the proof of Sophie
Germain's Theorem. Legendre does not mention the claims in the manuscript that
turn out not to be validly proved. If Legendre used Germain's Manuscript B as
his source for what he chose to publish as Germain's, then he vetted it and
extracted the parts that were correct.

\section{Even exponents\label{evenexponents}}

Another direction of Germain's is provided by three pages that we call
Manuscript C.\footnote{Yet one more manuscript, claiming to dispense with even
exponents by quite elementary means, is \cite[pp.~90v--90r]{ger9114}. It
contains a mistake that Germain went back to, crossed out, and corrected. But
she did not carry the corrected calculation forward, likely because it is then
obvious that it will not produce the desired result, so is not worth pursuing
further.} These pages contain highly polished statements with proof of two theorems.

The first theorem claims that the \textquotedblleft
near-Fermat\textquotedblright\ equation $2z^{m}=y^{m}+x^{m}$ (which amounts to
seeking three $m$-th powers in arithmetic progression) has no nontrivial
natural number solutions (i.e., other than $x=y=z$) for any even exponent
$m=2n$ with $n>1$. In fact Germain claims that her proof applies to an entire
family of similar equations in which the exponents are not always the same for
all variables simultaneously. Her proof begins with a parametric
characterization of integer solutions to the \textquotedblleft
near-Pythagorean\textquotedblright\ equation $2c^{2}=b^{2}+a^{2}$ (via
$c=z^{n}$, $b=y^{n}$, $a=x^{n}$), similar to the well-known parametric
characterization of Pythagorean triples (solutions to $c^{2}=b^{2}+a^{2}$)
used by Euler in his proof of Fermat's Last Theorem for exponent four \cite[p.
178]{lp}. The characterization of near-Pythagorean triples, stemming from a
long history of studying squares in arithmetic progression, would have been
well known at the time \cite[ch. XIV]{dickson}.

We will not analyze Germain's proof further here, nor pronounce judgement on
its correctness, except to say that it likely flounders in its fullest
generality near the beginning, as did the proof above of Theorem 3 in
Manuscript B, on another unjustified assumption of relative primality of two
expressions. However, this would still allow it to apply for \textquotedblleft
Case 1\textquotedblright, i.e., when $x$, $y$, $z$, are relatively prime to
$n$. Someone else may wish to pursue deciphering whether the entire proof is
valid in this case or not. There is a substantial history of research on the
near-Fermat equation $2z^{m}=y^{m}+x^{m}$. It was finally proven in 1997 by
Darmon and Merel \cite{darmon} to have no nontrivial solutions for $m>2$,
after partial results by Ribet \cite{ribet} and D\'{e}nes \cite{denes}, among
others. Much earlier, Euler had proved its impossibility for $m=4$
\cite{denes} \cite[ch. XXII]{dickson} \cite{ribet}, and then for $m=3$
\cite{denes} \cite[ch. XXI]{dickson}. So Germain's claim is now known to be
true, and it would be interesting to understand her method of proof well
enough to see if it is viable for Case 1.

Germain's second claim is to prove Fermat's Last Theorem for all even
exponents greater than two, i.e., for $z^{2n}=y^{2n}+x^{2n}$ with $n>1$, and
her proof relies directly on the previous theorem. It seems to us that this
proof too relies on the unsupported relative primality of two expressions, in
this case the two factors $z-y$ and $z^{n-1}+yz^{n-2}+\cdots+y^{n-2}z+y^{n-1}$
of $z^{n}-y^{n}$, under only the assumption that $x,$ $y,$ and $z$ are
pairwise relatively prime. It does seem to us that Germain's proof is fine,
though, for \textquotedblleft Case 1\textquotedblright\ (modulo appeal to the
previous theorem, of course), i.e., provided that $x$, $y$, $z$, are
relatively prime to $n$, in which case the two factors above will be
relatively prime. We note that it is under an almost identical hypothesis that
Terjanian proved Case 1 of Fermat's Last Theorem for even exponents in 1977
\cite[VI.4]{ribenboim-amateurs} \cite{terjanian}.

\section{Germain's approaches to Fermat's Last Theorem: pr\'{e}cis and
connections\label{precis}}

\begin{figure}
[pt]
\begin{center}
\includegraphics[
trim=0.000000in 0.000000in 0.004846in 0.000000in,
natheight=4.164100in,
natwidth=8.077300in,
height=2.6161in,
width=5.0462in
]%
{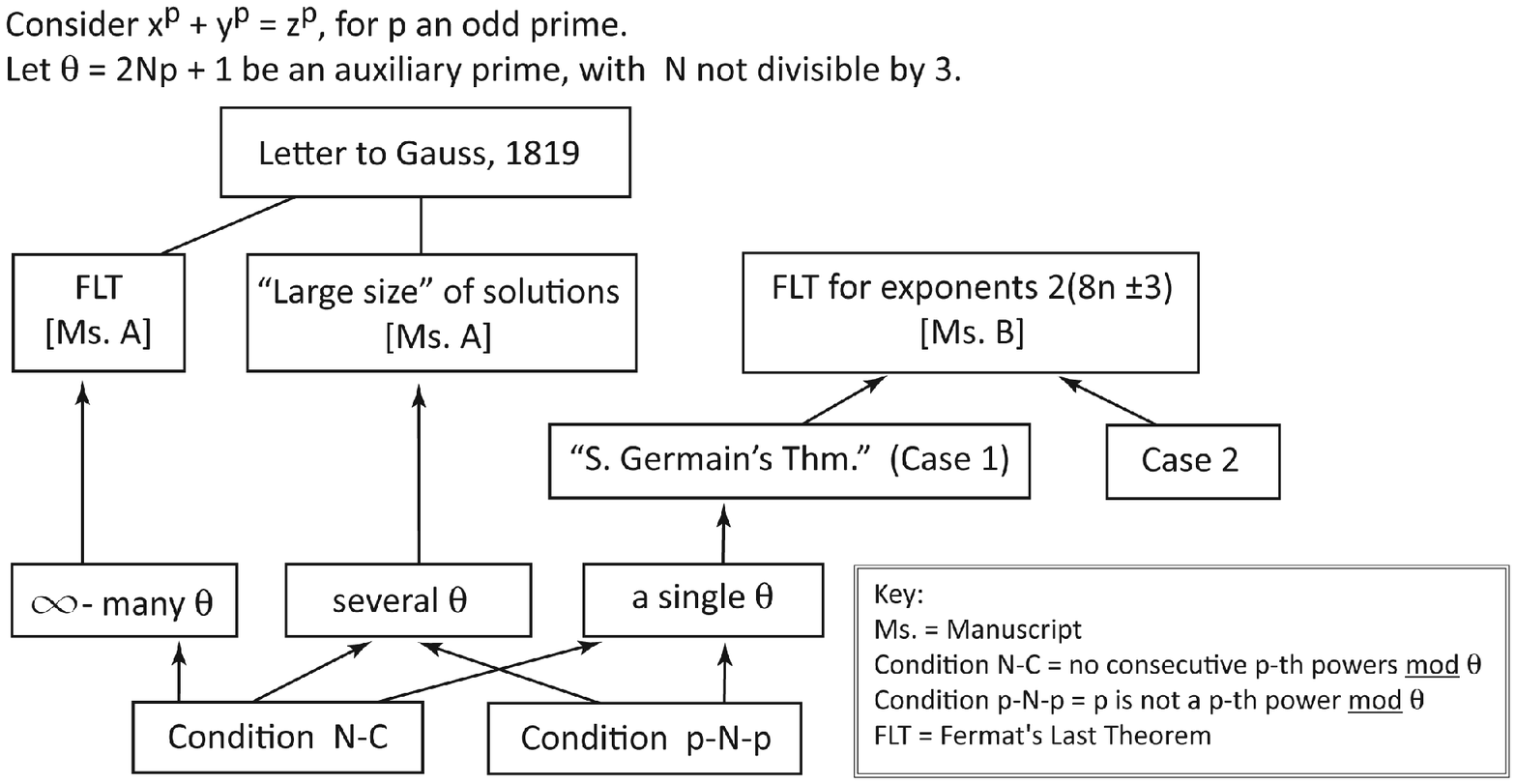}%
\caption{Conditions (hypotheses) for theorems}%
\label{figure-1}%
\end{center}
\end{figure}

\begin{figure}
[pt]
\begin{center}
\includegraphics[
natheight=3.602800in,
natwidth=5.752700in,
height=2.5425in,
width=4.0413in
]%
{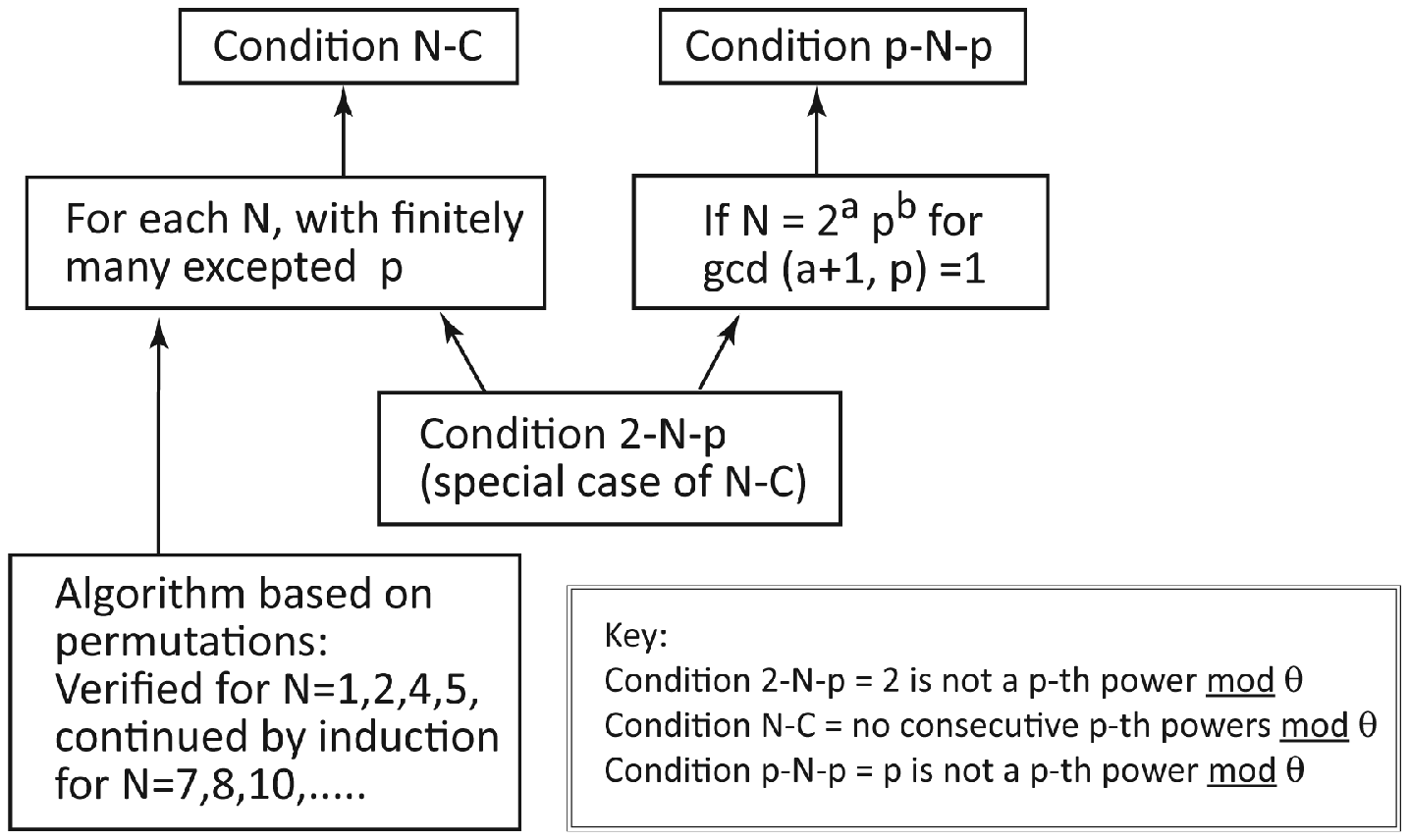}%
\caption{Algorithms and propositions for satisfying conditions}%
\label{figure-2}%
\end{center}
\end{figure}

Our analyses above of Sophie Germain's manuscripts have revealed a wealth of
important unevaluated work on Fermat's Last Theorem, calling for a
reassessment of her achievements and reputation. To prepare for our
reevaluation and conclusion, we first summarize (see Figures \ref{figure-1},
\ref{figure-2}) what we have discovered mathematically in these manuscripts,
and how it is related to other documentary evidence.%

%

\subsection{The grand plan to prove Fermat's Last Theorem}

In Manuscript A, Germain pioneers a grand plan for proving Fermat's Last
Theorem for any prime exponent $p>2$ based on satisfying a modular
non-consecutivity (N-C) condition for infinitely many auxiliary primes. She
develops an algorithm verifying the condition within certain ranges, and
outlines an induction on auxiliaries to carry her plan forward. Her techniques
for N-C verification are completely different from, but just as extensive as,
Legendre's, although his were for the purpose of proving Case 1, and were also
more ad hoc than hers. That Germain, as opposed to just Legendre, even had any
techniques for N-C verification, has been unknown to all subsequent
researchers who have labored for almost two centuries to extend N-C
verification for proving Case 1. Germain likely abandoned further efforts at
her grand plan after Legendre suggested to her that it would fail for $p=3$.
She sent him a proof confirming this, by showing that there are only finitely
many valid N-C auxiliaries.

Unlike Legendre's methods and terminology, Germain adopts Gauss's congruence
language and points of view from his \emph{Disquisitiones}, and thus her
techniques have in several respects a more group-theoretic flavor. Germain's
approach for verifying N-C was independently discovered by L. E. Dickson in
the twentieth century. He, or earlier researchers, could easily have obtained
a jump start on their own work by taking their cue from Germain's methods, had
they known of them. Recent researchers have again approached N-C by induction,
as did Germain.

\subsection{Large size of solutions and Sophie Germain's Theorem}

Also in Manuscript A, Germain writes a theorem and applications to force
extremely large minimal sizes for solutions to Fermat equations, based on
satisfying both the N-C and $p$-N-$p$ conditions. She later realized a flaw in
the proof, and attempted to repair it using her knowledge of quadratic
residues. The valid part of the proof yields what we call Sophie Germain's
Theorem, which then allows proof of Case 1 by satisfying the two conditions.

Germain's efforts to satisfy the $p$-N-$p$ condition are based on her
theoretical result showing that it will often follow from the $2$-N-$p$
condition, which she has already studied for N-C. This then makes it in
practice very easy to verify $p$-N-$p$, once again unlike Legendre. Germain's
result obtaining $p$-N-$p$ from $2$-N-$p$ was also independently discovered
much later, by Wendt, Dickson, and Vandiver in their efforts to prove Case 1.

\subsection{Exponents $2\left(  8n\pm3\right)  $ and Sophie Germain's Theorem}

In Manuscript B, Germain makes a very creditable attempt to prove Fermat's
Last Theorem for all exponents $2p$ where $p=8n\pm3$ is prime. Germain begins
with a proof of what we call Sophie Germain's Theorem, in order to argue for
Case 1. Manuscript B provides us with our best original source for the theorem
for which she is famous. Her subsequent argument for Case 2 boils down to
knowledge about biquadratic residues. This latter argument contains a flaw
related to relative primality. The manuscript fits well as a primary source
for what Legendre credited to Germain.

One could imagine that the appearance here of Sophie Germain's Theorem might
indicate an effort to recover what she could from the flawed Large Size
theorem in Manuscript A, but the details of the proof suggest otherwise, since
they betray the same misunderstanding as in Manuscript A before Germain wrote
its erratum.

\subsection{Even exponents}

In Manuscript C, Germain writes two theorems and their proofs to establish
Fermat's Last Theorem for all even exponents, by methods completely unlike
those in her other manuscripts. She plans to prove Fermat's Last Theorem by
showing first that a slightly different family of Diophantine equations has no
solutions. So she begins by claiming that the \textquotedblleft
near-Fermat\textquotedblright\ equations $2z^{2n}=y^{2n}+x^{2n}$ (and whole
families of related equations) have no nontrivial positive solutions for
$n>1$. This has only very recently been proven in the literature. Her proof
suffers from the same type of flaw for Case 2 as in Manuscript B, but may
otherwise be correct. Her proof of Fermat's Last Theorem for even exponents,
based on this \textquotedblleft near-Fermat result,\textquotedblright\ also
suffers from the Case 2 flaw, but otherwise appears to be correct.

\section{Reevaluation}

\subsection{Germain as strategist: theories and techniques}

We have seen that Germain focused on big, general theorems applicable to
infinitely many prime exponents in the Fermat equation, rather than simply
tackling single exponents as usually done by others. She developed general
theories and techniques quite multifaceted both in goal and methods. She did
not focus overly on examples or ad hoc solutions. And she also used to great
advantage the modern point of view on number theory espoused by Gauss. The
significance of Germain's theoretical techniques for verifying conditions N-C
and $p$-N-$p$ is indicated by their later rediscovery by others, and a recent
reapproach by mathematical induction. Moreover, her approach was more
systematic and theoretical than Legendre's pre-Gaussian and completely
different methods.

For almost two hundred years, Germain's broad, methodical attacks on Fermat's
Last Theorem have remained unread in her unpublished papers. And no one has
known that all the results published by Legendre verifying conditions N-C and
$p$-N-$p$, quoted and used extensively by others, are due but uncredited to
Germain, by more sophisticated and theoretical methods.

These features of Sophie Germain's work demonstrate that, contrary to what has
been thought by some, she was not a dabbler in number theory who happened to
light upon one significant theorem. In fact, what we call Sophie Germain's
Theorem is simply fallout from two much grander engagements in her papers,
fallout that we can retrospectively isolate, but which she did not. It is we
and Legendre, not Germain, who have created Sophie Germain's Theorem\ as an
entity. On the other hand, Legendre in this sense also performed a great
service to Germain and to future research, since he extracted from her work
and published the one fully proven major theorem of an enduring and broadly
applicable nature.

Germain's agenda was ambitious and bold. She tackled what we now know was one
of the hardest problems in mathematics. It is no surprise that her attempts
probably never actually proved Fermat's Last Theorem\ for even a single new
exponent, although she seems to have come close at times.

\subsection{Interpreting errors in the manuscripts}

Mathematicians often make errors in their work, usually winnowed out through
reactions to presentations, informal review by colleagues, or the publication
refereeing process. We have found that several of Germain's manuscripts on
Fermat's Last Theorem contain errors in her proofs. Let us examine these in
light of the unusual context within which we have found them.

First, we are short-circuiting normal publication processes by peeking at
Germain's private papers, works she chose never to submit for publication,
even had she shown them to anyone. Perhaps she knew of the errors we see, but
chose to keep these papers in a drawer for later revival via new ideas. We can
see explicitly that she later recognized one big error, in her Large Size of
Solutions proof, and wrote an erratum attempting remedy.

Second, let us consider the mathematical nature of the mistakes in her
manuscripts. In elasticity theory, where the holes in her societally forced
self-taught education were serious and difficult to remediate on her own
\cite[p. 40ff]{bucc}, Germain suffered from persistent conceptual difficulties
leading to repeated serious criticisms. By contrast, Germain was very
successful at self-education and independent work in number theory. She was
able to train herself well from the books of Legendre and Gauss, and she shows
careful work based on thorough understanding of Gauss's \emph{Disquisitiones
Arithmeticae}, despite its highly technical nature. The mistakes in her number
theory manuscripts do not stem from conceptual misunderstanding, but rather
are slips overlooking the necessity for relative primality in making certain
deductions, even though elsewhere she shows clear awareness of this necessity.
In particular, Germain's entire grand plan for proving Fermat's Last Theorem,
including algorithms for verifying Conditions N-C and $p$-N-$p$, is all very
sound. Even though Germain's mistakes were conceptually minor, they happen to
have left her big claims about large size and proving Fermat's Last Theorem
for various families of exponents unproven.

Further, we should ask what evaluation by peers Germain's manuscripts
received, that should have brought errors to her attention. Here we will
encounter more a puzzle than an answer.

\subsection{Review by others versus isolation}

\subsubsection{Germain's elasticity theory: praise and neglect}

There is already solid evidence \cite[Ch. 5--9]{bucc} that during Germain's
long process of working to solve the elasticity problem in mathematical
physics,\footnote{The Academy's elasticity prize competition was announced in
1809, twice extended, and Germain eventually received the award in 1816.
Thereafter she carried out efforts at personal, rather than institutional,
publication of her work on elasticity theory, stretching long into the 1820s.}
she received ever decreasing collegial review and honest critique of her work.
In fact, towards the end perhaps none.

Publicly praised as genius and marvel, she was increasingly ignored privately
and institutionally when it came to discourse about her elasticity work. There
is no evidence of any individual intentionally wishing her harm, and indeed
some tried personally to be quite supportive. But the existing system ensured
that she lacked early solid training or sufficiently detailed and constructive
critique that might have enabled her to be more successful in her research.
Germain labored continually under marginalizing handicaps of lack of access to
materials and to normal personal or institutional discourse, strictures that
male mathematicians did not experience \cite[Ch. 7--9]{bucc}. The evidence
suggests that Germain in effect worked in substantial isolation much of the time.

\subsubsection{Germain's interactions about Fermat's Last Theorem: the
evidence}

Given the social features dominating Germain's work in elasticity theory, what
was the balance between collegial interaction and isolation in her work?

Specifically, we will focus on what to make of the disparity between the
techniques of Germain and Legendre for their many identical results on the
Fermat problem. And we will ask what of Germain's work and results was seen by
Legendre, or anyone?

We have no actual published work by Germain on Fermat's Last Theorem. Even
though much of the research in her manuscripts would have been eminently
publishable, such as her theoretical means of verifying the N-C and $p$-N-$p$
conditions for applying Sophie Germain's Theorem to prove Case 1, it never
was. While we could speculate on reasons for this, it certainly means that it
did not receive any formal institutional review. Nor presumably could Germain
present her work to the Academy of Sciences, like her male contemporaries.

Despite having analyzed a wealth of mathematics in Germain's manuscripts, we
still have little to go on when considering her interactions with others. Her
manuscripts say nothing directly about outside influences, so we must infer
them from mathematical content.

Germain's 1819 letter to Gauss focused on the broad scope of her work on
Fermat's Last Theorem, but did not mention direct contact with others, and
apparently received no response from Gauss. Gauss had earlier made clear his
lack of interest in the Fermat problem, writing on March 21, 1816 to Olbers
\cite[p.~629]{olbers}: \textquotedblleft I am very much obliged for your news
concerning the [newly established] Paris prize. But I confess that Fermat's
theorem as an isolated proposition has very little interest for me, because I
could easily lay down a multitude of such propositions, which one could
neither prove nor dispose of.\textquotedblright\ This could by itself explain
why Germain did not receive a response from Gauss to her 1819 letter.

Thus the Fermat problem was in a very curious category. On the one hand, from
1816--1820 it was the subject of the French Academy's prize competition,
thereby perhaps greatly attracting Germain's interest. After all, with no
access to presenting her work at the Academy, her primary avenues for
dissemination and feedback were either traditional journal publication or the
Academy prize competition, which she had won in elasticity. On the other hand
the Fermat problem was considered marginal by Gauss and others, and topics
such as the investigation of higher reciprocity laws certainly involved
developing important concepts with much wider impact. So Germain's choice to
work mostly on Fermat's Last Theorem, while understandable, contributed to her
marginalization as well.

Regarding Germain's interaction with Legendre about her work on Fermat's Last
Theorem, we have two important pieces of evidence. First, while Legendre's
published footnote crediting Sophie Germain's Theorem to her is brief, we can
correlate it very precisely with content found in Germain's manuscripts.
Second, we have one critical piece of correspondence, Germain's letter to
Legendre confirming that her grand plan will not work. Starting from these we
will now draw some interesting conclusions.

\subsubsection{Legendre and Germain: A perplexing
record\label{Leg-Germ-interaction}}

Legendre's footnote and Germain's letter to him indicate that they had
mathematically significant contact about the Fermat problem, although we do
not know how frequently, or much about its nature. What then does our study of
her most polished manuscripts suggest?

First, it is a real surprise to have found from Manuscript A that Germain and
Legendre each had very extensive techniques for verifying Conditions N-C and
$p$-N-$p$, but that they are completely disjoint approaches, devoid of
mathematical overlap. Their methods were obviously developed completely
independently, hardly what one would expect from two mathematicians in close contact.

This phenomenon dovetails with a counterview about the effects of isolation
suggested to us by Paulo Ribenboim. If one works in isolation, one is not so
much influenced by others, so one has the advantage of originality, provided
one has fresh, good ideas. Clearly Germain had these, since we have seen that
she developed her own powerful theoretical techniques for verifying Conditions
N-C and $p$-N-$p$, not derived from anyone else's.

In contrast to Manuscript A, Legendre's crediting footnote details exactly the
results that are correct from Germain's Manuscript B, namely Sophie Germain's
Theorem and an additional technical result about the equations in its proof.
So while Manuscript B, along with her separate table of residues and
auxiliaries, is an extremely plausible source for Legendre's credit to her,
Germain's Manuscript A shows completely independent but parallel work left
invisible by Legendre's treatise.

So where does this leave Manuscript A? It contains Germain's grand plan, along
with all her methods and theoretical results for verifying N-C and $p$-N-$p$,
and her large size theorem. This seems like her most substantial work, and yet
we can find only a single speck of circumstantial evidence in Legendre's 1823
treatise suggesting that he might even be aware of the mathematics in
Germain's Manuscript A, despite her manuscript being placed by her letter to
Gauss at prior to 1819. But even this speck is perplexing and can be viewed in
opposing ways, as follows.

Recall from footnote \ref{legendre-to-1000} that Legendre, in his treatment of
large size of solutions, comments that for $p=5$ his data makes him
\textquotedblleft presume\textquotedblright\ that there are no auxiliary
primes larger than $101$ satisfying Condition N-C. This indicates that he was
at least interested in whether there are infinitely many auxiliaries, although
he does not mention why. Why would he even be interested in this issue, if it
weren't for interest in the grand plan? And why would he even imagine that
there might only be finitely many, unless he already had some evidence
supporting that, such as Germain's letter to him proving failure of the grand
plan for $p=3$? On the other hand, if he had her letter before writing his
1823 memoir, why did he not say something stronger for $p=5$, such as that he
knew that for $p=3$ there are only finitely many primes satisfying N-C,
supporting his presumption for $p=5$?

The only direct evidence we have that Legendre knew of Germain's grand plan is
her letter to him proving that it will not work for $p=3$. But even if
Germain's letter proving failure of the grand plan for $p=3$ occurred before
Legendre's 1823 treatise, so that the known failure was his reason for not
mentioning the plan in his treatise, why is Legendre mute about Germain
through the many pages of results identical to hers that he proves, by
completely different means, on Conditions N-C and $p$-N-$p$ for establishing
Case 1 and large size of solutions? Extensions of these results have been
important to future work ever since, but no one has known that these were
equally due to Germain, and by more powerful methods.

If Legendre had seen Manuscript A, he knew all about Germain's methods, and
could and should have credited her in the same way he did for what is in
Manuscript B. We must therefore at least consider, did Legendre, or anyone
else, ever see Manuscript A and so comprehend most of Germain's work, let
alone provide her with constructive feedback? It is reasonable to be
skeptical. Earlier correspondence with Legendre shows that, while he was a
great personal mentor to her initially during the elasticity competition, and
seems always to have been a friend and supporter, he withdrew somewhat from
mentorship in frustration as the competition progressed \cite[p.~63]{bucc}.
Did this withdrawal carry over somehow to contact about Fermat's Last Theorem?
Without finding more correspondence, we may never know whether Germain had
much extensive or intensive communication with anyone about her work on
Fermat's Last Theorem.

\subsubsection{The Fermat prize competition\label{prize-competition}}

There was one final possible avenue for review of Germain's work on the Fermat problem.

At the same session of the Academy of Sciences in 1816 at which Sophie Germain
was awarded the elasticity competition prize, a new competition was set, on
the Fermat problem. Extended in 1818, it was retired in 1820 with no award,
and Sophie Germain never made a submission \cite[p. 86]{bucc}. And yet,
together, our manuscript evidence and the 1819 date of her letter to Gauss
strongly suggest that she was working hard on the problem during the years of
the prize competition.

Why did she not submit a manuscript for this new prize, given the enormous
progress on the Fermat problem we have found in her manuscripts, and the
meticulous and comprehensive appearance of her work in Manuscript A, which
appears prepared for public consumption? Was Germain's reluctance due to
previous frustrating experiences from her multiple submissions for the
elasticity prize through its two extensions---a process that often lacked
helpful critiques or suggested directions for improvement \cite[Ch.
5--9]{bucc}? Or, having been particularly criticized for incompleteness during
the elasticity prize competition, did she simply know she had not definitely
proved Fermat's Last Theorem in full, and hence felt she had nothing
sufficient to submit?

\subsection{Amateur or professional?}

Goldstein \cite{goldstein-amateur} analyzes the transformation of number
theory from the domain of the amateur to that of the professional during the
17th to 19th centuries. By Germain's time this transformation had shifted
number theory mostly to the professional world, and to be successful Germain
needed to interact and even compete with degreed professionals at
institutions. Was she herself an amateur or a professional?

Germain had many of the characteristics of a professional, attained through
highly unusual, in fact audacious, personal initiatives injecting herself into
a professional world that institutionally kept her, as a woman (and therefore
by definition uncertified), at arm's length. Her initiatives would hardly be
dreamt of by anyone even today. She attained some informal university
education first through impersonation of LeBlanc, a student at the \'{E}cole
Polytechnique, an institution that would not admit women, leading to
mathematicians like Lagrange and Legendre serving as her personal mentors. She
devoured much professional mathematical literature in multiple disciplines, to
which however she presumably had only what access she could obtain privately.
And she initiated an also impersonated correspondence with Gauss. Germain
appears to have devoted her adult life almost entirely to mathematical
research, having no paid employment, spouse, or children. She competed against
professional mathematicians for the Academy prize on elasticity, she achieved
some professional journal publications, and she self-published her elasticity
prize research when the Academy would not.

On the other hand, Germain had some of the characteristics of amateurs typical
of earlier periods, such as great reliance on personal contact and letters.
Most importantly, she was not employed as a professional mathematician. And
after her death no institution took responsibility for her papers or their
publication, one substantial reason why much of her extensive work has
remained unknown. However, it seems that all this was ultimately due precisely
to her being a woman, with professional positions closed to her. One could say
that Germain was relegated to something of the role of an amateur by a world
of professionals and institutions that largely excluded her because of her
sex, a world to which she aspired and for which she would have otherwise been
perfectly qualified.

\section{Conclusion}

The impression to date, the main thesis of \cite{bucc}, has been that Germain
could have accomplished so much more had she enjoyed the normal access to
education, collegial interaction and review, professional institutions, and
publication accorded to male mathematicians. Our study of her manuscripts and
letters bolsters this perspective.

The evidence from Germain's manuscripts, and comparison of her work with that
of Legendre and later researchers, displays bold, sophisticated, multifaceted,
independent work on Fermat's Last Theorem, much more extensive than the single
result, named Sophie Germain's Theorem, that we have had from Legendre's
published crediting footnote. It corroborates the isolation within which she
worked, and suggests that much of this impressive work may never have been
seen by others. We see that Germain was clearly a strategist, who
single-handedly created and pushed full-fledged programs towards Fermat's Last
Theorem, and developed powerful theoretical techniques for carrying these out,
such as her methods for verifying Conditions N-C and $p$-N-$p$.

We are reminded again of her letter to Gauss: \textquotedblleft I will give
you a sense of my absorption with this area of research by admitting to you
that even without any hope of success, I still prefer it to other work which
might interest me while I think about it, and which is sure to yield
results.\textquotedblright\footnote{(Letter to Gauss, p. 2)}\ Sophie Germain
was a much more impressive number theorist than anyone has ever known.

\section*{Acknowledgements}

We owe heartfelt thanks to many people who have helped us tremendously with
this project since 1993: Evelyne Barbin, H\'{e}l\`{e}ne Barcelo, Louis
Bucciarelli, Keith Dennis, Mai Gehrke, Tiziana Giorgi, Catherine Goldstein,
Maria Christina Mariani, Pat Penfield, Donato Pineider, Paulo Ribenboim, and
Ed Sandifer, along with Marta Gori of the Biblioteca Moreniana, as well as the
Biblioth\`{e}que Nationale, New York Public Library, Nieders\"{a}chsische
Staats- und Universit\"{a}tsbibliothek G\"{o}ttingen, and the Interlibrary
Loan staff of New Mexico State University. We also thank the referees and
editors for thoughtful, provocative, and tremendously helpful critique and suggestions.

We thank the Biblioth\`{e}que Nationale de France, the Nieders\"{a}chsische
Staats- und Universit\"{a}tsbibliothek G\"{o}ttingen, and the New York Public
Library, for permission to reproduce from manuscripts and letters.

\end{document}